\documentclass[reqno,10pt,centertags,draft]{amsart}
\usepackage{amsmath,amsthm,amscd,amssymb,latexsym,upref}
\date{\today}

\newtheorem{theorem}{Theorem}[section]
\newtheorem{proposition}[theorem]{Proposition}
\newtheorem{lemma}[theorem]{Lemma}
\newtheorem{corollary}[theorem]{Corollary}

 \theoremstyle{definition}
 \newtheorem{definition}[theorem]{Definition}
\newtheorem{remark}[theorem]{Remark}

\newcommand{\ZZ}{{\mathbb Z}}
\newcommand{\RR}{{\mathbb R}}
\newcommand{\CC}{{\mathbb C}}
\newcommand{\TT}{{\mathbb T}}
\newcommand{\EE}{{\mathbb E}}
\newcommand{\DD}{{\mathbb D}}

\newcommand{\fA}{{\mathfrak A}}
\newcommand{\fL}{{\mathfrak L}}
\newcommand{\bphi}{{\boldsymbol{\phi}}}
\newcommand{\fR}{{\mathfrak R}}

\newcommand{\cL}{{\mathcal L}}

\newcommand{\hJ}{\widehat{J}}
\newcommand{\hS}{\widehat{S}}

\newcommand{\nLp}{{L^2_{d\chi}(l^2(\ZZ_+))}}
\newcommand{\Julia}{{\rm Jul}}
\newcommand{\Jac}{{\rm Jac}}
\renewcommand{\div}{{\rm div}}
\newcommand{\Drb}{{\rm Drb}}

\title[  Finite difference operators with a finite--band spectrum]
{Finite difference operators with a finite--band spectrum}

\author[]{ F. Peherstorfer,  and P. Yuditskii}

\thanks{
Partially supported by   the Austrian Founds FWF, project number:
P16390--N04 and {\it Marie Curie International Fellowship} within
the 6-th European Community Framework Programme,
 Contract MIF1-CT-2005-006966}
\begin{document}

\maketitle


\begin{abstract}
We study  the correspondence between almost periodic difference
operators and algebraic curves (spectral surfaces). An especial role
plays the parametrization of the spectral curves in terms of, so
called, branching divisors.  The multiplication operator by the
covering map with respect to the natural basis in the Hardy space on
the surface is the $2d+1$--diagonal matrix; the $d$--root of the
product of the Green functions (counting their multiplicities) with
respect to all infinite points on the surface is the symbol of the
shift operator. We demonstrate an application of our general
construction to a particular covering, which generate  widely
discussed almost periodic  CMV matrices. We discuss an important
theme: covering of one spectral surface by another one and related
to this operation transformations on the set of multidiagonal
operators (so called Renormalization Equations). We proof several
new results dealing with Renormalization Equations for periodic
Jacobi matrices (polynomial coverings) and the case of  a rational
double covering.
\end{abstract}

\section{Introduction}

\subsection{Ergodic finite difference operators and associated
Riemann surfaces}

The standard (three--diagonal) finite--band Jacobi matrices
\cite{DKN, Tes} can be defined as almost periodic or even ergodic
Jacobi matrices with absolutely continuous spectrum that consists of
a finite system of intervals. We wish to find a natural extension of
this class of finite difference operators to the multi-diagonal
case.

The correspondence between periodic difference operators and
algebraic curves (spectral surfaces) was discussed in detail in
\cite{MM}. Let us generalize this construction to show (at least on
a speculative level) how almost periodic or ergodic operators give
rise to a corresponding spectral surface.

Recall the definition of an ergodic operator \cite{Cy, PF}, see also
\cite{VYu}. Let $(\Omega,\frak A, d\chi)$ be a separable probability
space and let $\mathcal T:\Omega\to\Omega$ be an invertible ergodic
transformation, i.e., $\mathcal T$ is measurable, it preserves
$d\chi$, and every measurable $\mathcal T$--invariant set has
measure $0$ or $1$. Let $\{q^{(k)}\}_{k=0}^d$ be functions from
$L^\infty_{d\chi}$, with $q^{(d)}$  positive--valued and $q^{(0)}$
real--valued. Note that in the periodic case $\Omega=\ZZ/N\ZZ$,
where $N$ is the period and $\mathcal T\{n\}=\{n+1\}$, $\{n\}\in
\ZZ/N\ZZ$.

Then with almost every $\omega\in\Omega$ we associate a
self--adjoint $2d+1$--diagonal operator $J(\omega)$ as follows:

\begin{equation}\label{general}
(J(\omega)x)_n= \sum_{k=-d}^{d}
\overline{q^{(k)}_{n}(\omega)}x_{n+k}, \quad
          x=\{x_n\}_{n=-\infty}^{\infty} \in l^2(\ZZ),
\end{equation}
where
  $q_n^{(k)}(\omega):=q^{(k)}(\mathcal T^n\omega)$
and $q^{(-k)}(\omega):= \overline{q^{(k)}
(\mathcal T^{-k}\omega)}$.

Note that the structure of $J(\omega)$ is described by the
following identity
\begin{equation}\label{com}
J(\omega)S=S J(\mathcal T\omega),
\end{equation}
where $S$ is the shift operator in $l^2(\ZZ)$. The last relation
indicates strongly that one can associate with the family of
matrices $\{J(\omega)\}_{\omega\in\Omega}$ a natural pair of
commuting operators (in the periodic case one just uses the fact
that $J$ and $S^N$ commute).

Namely, let $L^2_{d\chi}(l^2(\ZZ))$ be the space of
$l^2(\ZZ)$--valued vector functions, $x(\omega)\in l^2(\ZZ)$, with
the norm
\begin{equation*}
||x||^2=\int_\Omega||x(\omega)||^2\,d\chi.
\end{equation*}
Define
\begin{equation*}
(\hJ x)(\omega)=J(\omega)x(\omega), \quad (\hS x)(\omega)=S
x(\mathcal T\omega), \quad x\in L^2_{d\chi}(l^2(\ZZ)).
\end{equation*}
Then \eqref{com} implies
\begin{equation*}
(\hJ\hS x)(\omega)= J(\omega)S x(\mathcal T\omega)= S
J(\mathcal T\omega)
x(\mathcal T\omega)= (\hS\hJ x)(\omega).
\end{equation*}

Further, $\hS$ is a unitary operator and
$\hJ$ is self-adjoint. The space
$L^2_{d\chi}(l^2(\ZZ_+))$
is an invariant subspace for $\hS$. It is not
invariant with respect to $\hJ$ but it is invariant with respect to
$\hJ\hS^d$.
Put
\begin{equation*}
\hS_+=\hS| L^2_{d\chi}(l^2(\ZZ_+)),
\quad
(\hJ\hS^d)_+=\hJ\hS^d| L^2_{d\chi}(l^2(\ZZ_+)).
\end{equation*}


\begin{definition}[local functional model]
We say that a pair of commuting operators $A_1:H\to H$ and $A_2:H\to
H$ has a (local) functional model if there is a unitary embedding
$i:H\to H_O$ in a space $H_O$ of functions $F(\zeta)$ holomorphic in
some domain $O$ with a reproducing kernel ($F\mapsto F(\zeta_0)$,
$\zeta_0\in O$, is a bounded functional in $H_O$) such that
operators $i_\star A_1$ and $i_\star A_2$ become a pair of operators
of multiplication by holomorphic functions, say
\begin{equation*}
A_1 x\mapsto a_1(\zeta) F(\zeta), \quad A_2 x\mapsto a_2(\zeta)
F(\zeta).
\end{equation*}
\end{definition}

Existence of a local functional model implies a number of quite
strong consequences. In what follows $a(\zeta)$  and $b(\zeta)$
denote the functions (symbols) related to the operators
$(\hJ\hS^d)_+$
 and $\hS_+$. Let $k_\zeta$ be the reproducing
kernel in $H_O$ and let $\hat k_\zeta$ be its preimage $i^{-1}
k_\zeta$ in $L^2_{d\chi}(l^2(\ZZ_+))$. Then
\begin{equation}\label{3}
\langle \hS_+^* \hat k_\zeta, x\rangle
=\langle \hat k_\zeta, \hS_+ x\rangle=
\langle k_\zeta, b F\rangle.
\end{equation}
By the reproducing property
\begin{equation*}
\langle k_\zeta, b F\rangle=
\overline{b(\zeta) F(\zeta)}=
\langle \overline{b(\zeta)}k_\zeta,  F\rangle.
\end{equation*}

Hence,
\begin{equation*}
\langle \hS_+^* \hat k_\zeta, x\rangle=
\langle \overline{b(\zeta)}\hat k_\zeta,
x\rangle.
\end{equation*}
That is $\hat k_\zeta$ is an eigenvector of $\hS^*_+$ with the
eigenvalue $\overline {b(\zeta)}$. In the same way, $\hat k_\zeta$
is an eigenvector of $(\hJ\hS)^*_+$ with the eigenvalue
$\overline{a(\zeta)}$.

Thus, if a functional model exists then
the spectral problem
\begin{equation}\label{4}
\left\{\begin{matrix}
  \hS_+^* \hat k_\zeta &=
\overline{b(\zeta)}\hat k_\zeta\\
(\hJ\hS^d)_+^* \hat k_\zeta &=
\overline{a(\zeta)}\hat k_\zeta
\end{matrix},\right.
\end{equation}
has a solution $\hat k_\zeta $
antiholomorphic in $\zeta$.
 Moreover, linear combinations of all $\hat
k_\zeta$ are dense in $L^2_{d\chi}(l^2(\ZZ_+))$. Vice versa, if
\eqref{4} has a solution of such kind  then we define
\begin{equation*}
F(\zeta):=\langle x, \hat k_\zeta\rangle,
\quad ||F||^2:=||x||^2.
\end{equation*}
This provides a local functional
model for the pair  $\hS_+$,
  $(\hJ\hS)_+$.

The following proposition is evident.
\begin{proposition} Let $\mathcal U:L^2_{d\chi}\to
L^2_{d\chi}$ be the unitary operator associated with the ergodic
transformation $\mathcal T$: $(\mathcal Uc)(\omega)=
c(\mathcal T\omega)$,
$c\in L^2_{d\chi}$.
We denote by the
same letter $q$ both a function $q\in L^\infty_{d\chi}$ and the
 multiplication
operator $q$ (e.g., $(q c)(\omega):= q (\omega) c(\omega)$). Problem
\eqref{4} is equivalent to the following spectral problem
\begin{equation}\label{5}
\left\{
\sum_{k=-d}^d \mathcal U^k
\overline{q^{(k)}b^k(\zeta)}\right\}
c_\zeta=\overline{z(\zeta)}c_\zeta,
\end{equation}
where $z(\zeta):=a(\zeta)/b^d(\zeta)$ and $c_\zeta$ is an
anti--holomorphic $L^2_{d\chi}$--valued vector function. Moreover
$\{c_\zeta\}$ is complete in $L^2_{d\chi}$ if and only if $\{\hat
k_\zeta\}$ is complete in $L^2_{d\chi}(l^2(\ZZ_+))$.
\end{proposition}

We may hope to glue the local functional realizations to a {\it
global} functional model on a Riemann surface $X_0=\DD/\Gamma_0$
formed by functions $(z,b)$. Of course, this model does not
necessarily exist (even existence of a local model requires some
additional assumptions on the ergodic map and the coefficients
functions).

But, in particular, in the periodic case, when
$L^2_{d\chi}=\CC^N$,
we have
$$
\mathcal U\begin{bmatrix} c_0\\ \vdots \\c_{N-1}
\end{bmatrix}=
\begin{bmatrix} c_1\\ \vdots \\c_{0}
\end{bmatrix},
\quad
\mathcal U=\begin{bmatrix} 0&1& & \\ & \ddots &\ddots& \\ & & \ddots &
1\\ 1& & & 0
\end{bmatrix},
$$
where $c_n:=c(\{n\})$ and the $q^{(k)}$'s became the diagonal
matrices. Thus \eqref{5} means that the $N\times N$ matrix has a
nontrivial annihilating vector and we arrive at the curve $X_0$ in
the form:
\begin{equation}\label{5.06}
\det\left[
\sum_{k=-d}^d \mathcal U^k
\overline{q^{(k)}b^k}-\overline{z} I\right]
=0
\end{equation}
and the restriction $|b|<1$.

The surface $X_0$ is in generic case of infinite genus. However we
can reduce it because $X_0$ possesses a family of automorphisms. Let
$e_{\{\gamma\}}$ be an eigenvector of $\mathcal U$ with an
eigenvalue $\bar\mu_{\{\gamma\}}$. The systems of eigenfunctions and
eigenvalues form both Abelian groups with respect to multiplication.
Using \eqref{5} we get immediately that $\{\gamma\}:(z,b)\mapsto
(z,\mu_{\{\gamma\}}b)$ is an automorphism of $X_0$. Taking a
quotient of $X_0$ with respect to these automorphisms we obtain a
much smaller surface $X=\DD/\Gamma$, $\Gamma=\{\gamma\Gamma_0\}$. In
periodic case this means that \eqref{5.06} is actually a polynomial
expression in $z$, $\lambda$ and $\lambda^{-1}$ where $\lambda:=b^N$
(since it is invariant with respect to the substitutions $b\mapsto
e^{2\pi i\frac{ k}{N}}b$, $0<k<N$), see \cite{MM}, see also Sect. 6.

  Note that $z$ is still a function on $X$ but
$b$ becomes a character automorphic function. Finally, using $z$ we
may glue the boundary of $X$, remove punched points (where
$z=\infty$) and get in this way a compact
 Riemann surface $X_c$ , such that
 $$X=(X_c\setminus\{P:z(P)=\infty\})\setminus E.$$

 The simplest
assumption is that the boundary $E$ is a finite system of cuts on
$X_c$. In this case we get that the  triple $\{X_c,z,E\}$
characterizes the spectrum  of a finite difference operator.  We
came to this triple basically due to heuristic arguments, but the
opposite direction is already a solid mathematical fact:
 every triple of this kind gives rise to a family of
ergodic finite difference operators \cite[Sect. 6]{MM}. In Sect. 2
and 3 we give details on constructions of such operators using the
theory  of Hardy spaces on Riemann surfaces.

Now we would like to note that one  can describe all triples of a
given type up to a natural equivalence relation. This natural
parametrization of the triples is one of the main point of the
current paper. It allows us to put consideration from the pure
algebraic points of view to a wider setting and to recruit very much
analytic tools.

\subsection{Parametrization of the spectral curves in terms of
branching divisors} Probably the best known result in  spectral
theory  of the nature we want to discuss  deals with the description
of the spectrum of periodic Jacobi matrices. Such a set $E$ should
have a form of an inverse polynomial image  $E= T^{-1}([-1,1])$. The
polynomial $T$ should have all critical points real $\{c_k:
T'(c_k)=0\}\subset \RR$, moreover all critical values $t_k:=T(c_k)$
should have modulus not less than $1$ and their signs should
alternate, i.e.: $|T(c_k)|\ge 1$ and
$$
T(c_{k-1})T(c_k)<0\quad\text{ for}\ \
\ldots < c_{k-1}< c_k<
\ldots.$$
The claim is that the system of numbers $\{t_k:|t_k|\ge 1, \
t_{k-1}t_k<0\},$ determines a polynomial with the prescribed and
ordered critical values $t_k$ uniquely,
modulo a change of the independent
variable $z\mapsto az+b,\; a>0, b\in\RR$. The proof uses a special
representation of $T$:
$$
T(z)=\pm\cos\phi(z),
$$
where $\phi$ is the conformal map of the upper half plane onto the half
strip with a system of cuts:
$$
\Pi=\{w=u+iv: 0\le v\le \pi d\}\setminus
\cup_{k}\{ v=\pi k, \ u\le h_k\},
$$
where $d=\deg T$ and $\cosh(h_k):=|t_k|$, $\phi(\infty)=\infty$.

We mention a more general theorem of Maclane \cite{Mc} and Vinberg
\cite{Vin} on the existence and uniqueness of real polynomials
(actually, and entire functions) with prescribed (ordered!)
sequences of critical values. In the case of polynomials, this
theorem says that there is a one to one correspondence between
finite ``up-down'' real sequences $$\ldots \leq t_{k-1}\geq t_k\leq
t_{k+1}\geq \ldots,$$ and real polynomials whose all critical points
are real, also, modulo a change of the independent variable
$z\mapsto az+b,\; a>0, b\in\RR$. The MacLane--Vinberg theorem is
based on an explicit description of the Riemann surfaces spread over
the plane of the inverse functions $T^{-1}$.

To be more precise in  the general case we start with the following
\begin{definition}\label{def21}
We say that two triples $({X_c}_1,z_1,E_1)$
and $({X_c}_2,z_2,E_2)$ are equivalent
if there exists a holomorphic homeomorphism
$h:{X_c}_1\to {X_c}_2$ such that
$z_1=h^*(z_2)$ and $E_2=h(E_1)$.
\end{definition}

Note, that for any triple $(X_c,z,E)$ the holomorphic function
$z:X_c\to\overline\CC$, $\overline\CC=\CC\cup\{\infty\}$, is a
ramified covering of $\overline\CC$. The fact is that it is possible
and  very convenient to describe equivalence classes of ramified
coverings in terms of branching divisor. Namely, a point $P\in X_c$
such that $\frac{dz}{d \zeta}|_P=0$ where $\zeta$ is a local
holomorphic coordinate in a neighborhood of $P$ is called
\emph{ramification point (or, critical point)}. Its image $z(P)$ is
called \emph{a branching point (or, critical value)}. The set of all
branching points of function $z$ form a branching divisor of $z$.
Note that infinity also can be a branching point. Since in our
consideration it plays an exclusive role it is convenient  in what
follows to denote by $\mathcal Z:= \{z_i\}_{i=1}^N$ all branching
points from the {\it finite part} of the complex plane $\CC$.

Clearly, branching divisors of
equivalent functions are the same. Moreover,
the compact holomorphic curve $X_c$ is also
uniquely determined by the branching divisor
and some additional ramification data (of
combinatorial type). Namely, assume that $z$
has degree $d$.
   Let $w_0\in \CC$ be a non
branching point. Fix a system of non-intersecting paths
$\gamma=\{\gamma_i\}_{i=1,\dots,N}$. The $i$th path $\gamma_i$
connects $w_0$ and $z_i$. We want to construct a system of loops
$l_i\subset \CC$. To construct $l_i$ we start from $w_0$ and follow
first $\gamma_i$ almost to $z_i$, then encircle $z_i$
counterclockwise along a small circle and finally go back to $w_0$
along $-\gamma_i$. Using $l_i$ we associate with each branching
point an element of the permutation group $\sigma_i\in \Sigma_d$.

The point $w_0$ has exactly $d$ preimages. Let us label them by
integers $\{1,\dots,d\}$. Let us follow the loop $l_i$ and lift this
loop to $X_c$ starting from each of the preimages of $w_0$. The
monodromy along path $l_i$ gives us a permutation
$\sigma_i\in\Sigma_d$ of preimages. We add to this system of loops
one more $l_{\infty}$, related to infinity (starting from $w_0$ we
go sufficiently far then make a big circle in the clockwise
direction and go back to $w_0$). This loop gives us one more
permutation $\sigma_{\infty}$.

Note that the product
$\sigma_1\cdot\dots\ \cdot\sigma_N$ times $\sigma_{\infty}$ is the
identity operator.
 Therefore, the function $z$, including its behavior at infinity,
determines $N$ branching points and $N$ permutations. These
permutations are not uniquely defined, they depend on the labeling
of preimages of $w_0$. Therefore, they are determined, up to a
conjugacy by the elements of $\Sigma_d$.

Given a set of branching points
$\mathcal Z(z)=\{z_i\}_{i=1,\dots,N}\subset\CC $ and a system of
permutations $\sigma(z)=(\sigma_1,\dots,\sigma_N)\in
\Sigma_d\times\dots\times \Sigma_d/ \Sigma_d$,  where the last quotient
is taken
 with respect to diagonal conjugation,
we can restore by Riemann theorem
the surface $X_c$ and the function $z$. Actually there is one more
topological condition: the surface should be connected. Throughout
the paper we assume that the system of permutations guaranteed this
condition to be hold.

Hence, the triple $(X_c,z,E)$ is equivalent
to the triple $(\mathcal Z,\sigma_\gamma,E)$.
We use this  triple as free parameters determining the spectral
surface of a $2d+1$-diagonal matrix.

\smallskip

\noindent {\bf Summary.} Comparably with the case of Jacobi
matrices, where we have only system of cuts (spectral intervals) in
the complex plane, in the multidiagonal case we have a new
additional system of parameters. We have to fix in $\CC$ a system of
critical points $\mathcal Z$, and associate to them  a system of
permutation $\sigma$, which actually depends on the base point
$w_0\in \CC\setminus \mathcal Z$ and the system of paths $\gamma$.
They define a  Riemann surfaces $X_c$ and a covering $z$. Then on
the set $z^{-1}(\RR)\subset X_c$ we chose a system of cuts $E$, and
thus $(X_c,z,E)$ is restored up to the equivalence \ref{def21}.

\begin{remark} As it was mentioned by the branching divisors language
one can extend   consideration from the pure algebraic level. In
particular, we are very interested in  infinite dimensional
generalizations: the point is that starting from the $5$--diagonal
case we have the branching divisor $\mathcal Z$ as a completely new
system of parameters characterizing the spectral surface. What is
its influence on the properties of the corresponding surface?
Consider, indeed, the simplest case of $5$--diagonal matrices. Then
we have just to specify the point set $\mathcal Z$ (all permutations
are of the form $\sigma_j=\left(\begin{matrix} 1&
2\\2&1\end{matrix}\right)$). What is the speed of accumulation of an
{\it infinite} system of points $\mathcal Z$ to, say,  {\it a
finite} system of cuts $E$ so that $X_c\setminus E$ is of Widom
type, or of Widom type with the Direct Cauchy Theorem  (for the
definition of these type of surfaces see \cite{Has}, see also
\cite{SYU1} for their role in the spectral theory)? We provide here
the following {\it example}. Let $W(z)$ be the infinite Blaschke
product in the upper half plane with zeros at $\{z_k\}_{k\ge 0}$.
Define
$$
\fR=\{P=(z,w):w^2=W,\quad \Im z\ge 0\}.
$$
Then $\sqrt{\frac{z-z_0}{z-\bar z_0}}$ is the Green function in
$\fR$. Note that $C_k=(z_k,0)$, $k\ge 1$ are its critical  points,
and therefore the Widom function is
$$
\Delta=\sqrt{W\frac{z-\bar z_0}{z- z_0}}.
$$
Note that the Carleson condition for this function on $\fR$ is the
standard Carleson condition for $\{z_k\}_{k\ge 1}$ in the upper half
plane. Thus  the Blaschke  and Carleson conditions on  the zero set
in the half plane guarantee that $\fR$ is of Widom type and of Widom
type with the Direct Cauchy Theorem respectively.
\end{remark}

\subsection{Structure of the paper and main results}
In Sect. 2, having the spectral surface fixed, we define a system
 of Hardy spaces on it (natural counterparts of the Hardy space in the
unit disk). In that section we assume that there is only one
"infinity" $P_0\in X_c\setminus E$, $z(P_0)=\infty$, that is, the
product of permutations $\sigma_1\cdot\dots\, \cdot\sigma_N$ is a
cycle. There is an intrinsic basis in the Hardy space: each next
basis elements has  at $P_0$ a zero of bigger and bigger
multiplicity. Of course, this is a counterpart of the standard basis
system $\{\zeta^n\}_{n\ge 0}$ in the standard $H^2$. Extending this
basis to the negative integers (the system extends in the direction
of functions having a pole at $P_0$ with growing multiplicity) one
gets a basis in the whole $L^2$. Finally the multiplication operator
by the covering map with respect to this basis is the
$2d+1$--diagonal matrix (to this end it is important to note that
$z$ has pole of multiplicity $d$ at $P_0$). The complex Green
function with zero at $P_0$ is playing the role of the symbol of the
shift operator. Then we study the question of uniqueness of such a
model for an ergodic operator (see Theorem \ref{model} and the
example right after it).

The general case (several "infinities") is considered in Sect. 3.
The $d$--root of the product of the Green functions (counting their
multiplicities) with respect to all infinite points on the surface
becomes the symbol of the shift.  We have the ordering of the
infinities as one more parameter, defining the basis system and the
corresponding multidiagonal matrix.

In Sect. 4 we demonstrate, how an application of our general
construction to a particular covering, generate the well known and
now widely discussed almost periodic  CMV matrices.

Starting from Sect. 5 we discuss an important theme: covering of one
spectral surface by another one and related to this operation
transformations on the set of multidiagonal operators (so called
Renormalization Equations). Let $\pi: Y_c\to X_c$ be a $d$--sheeted
covering. For a system of cuts $E$ define $F=\pi^{-1}(E)$. Then we
have $\pi: Y_c\setminus F\to X_c\setminus E$. The study is based on
the relation between the Hardy spaces on these Riemann surfaces.

The Renormalization Equations generated by polynomial coverings have
play\-ed an important role in studying of almost periodic Jacobi
matrices with a singular continuous spectrum, \cite{B}, see also
\cite{PVYU}, \cite{VY}. They act in the most natural way on periodic
Jacobi matrices, see Sect. 6.

In Sect. 7 we proof several new results dealing with Renormalization
Equations for periodic Jacobi matrices: describe the complete set of
their solutions; show their relation with the Ruelle operators.
Finally, we give a possible generalization of the constructions from
Sect. 6 for a wider class of almost periodic Jacobi matrices with a
singular continuous spectrum. In particular, we prove the Lipschitz
property of the Darboux transform.

Having in mind importance of the paper \cite{BBM}, where the
Renormalization Equation  generated by "just" {\it quadratic
polynomial} was used, we investigate in Sect. 8 the case of {\it
rational double covering} $\pi(v)=\tau v-\frac{\tau-1}{v}$,
$\tau>1$. As usual, the renormalization procedure is simpler to
formulate for operators acting on the (integer) half axis.
\begin{definition}\label{22def}
Let $A$ be a self--adjoint operator  acting in $l_+^2=l^2(\ZZ_+)$
with a cyclic vector $|0\rangle$ and the spectrum on $[-1,1]$. We
define its transform $\pi^*(A)$ in the following steps. First we
define the upper triangular matrix $\Phi$ (with positive diagonal
entries) by the condition
\begin{equation}
    A^2+4\tau(\tau-1)=\Phi^* \Phi.
\end{equation}
Introduce $A_*:=\Phi A\Phi^{-1}$ and define the operator
\begin{equation}
\left[\begin{array}{cc}
    A&\Phi^*\\ \Phi & A_*
    \end{array}\right],
\end{equation}
acting in $l^2_+\oplus l^2_+$. Finally, using the unitary operator
$U: l^2_+\to l^2_+\oplus l^2_+$, such that
\begin{equation}
    U|2k\rangle= |k\rangle\oplus 0,\quad
        U|2k+1\rangle= 0\oplus |k\rangle
\end{equation}
we construct
\begin{equation}
\pi^*(A):= \frac 1 {2\tau}U^* \left[\begin{array}{cc}
    A&\Phi^*\\ \Phi & A_*
    \end{array}\right]U: l^2_+\to l^2_+.
\end{equation}
\end{definition}
We have a theorem (see Theorem \ref{t27n}) on the weak convergence
of the iterative procedure $A_{n+1}=\pi^*(A_n)$ to an  operator with
a simple singular continuous spectrum supported on the Julia set of
the given expanding mapping and pose here a question on the
contractivity of the renorm operator, at least for  big values of
$\tau$. The main general conjecture deals with contractivity of all
renormalizations, generated by a covering with sufficiently big
critical values.

\section{The Global Functional Model (single infinity case).
Uniqueness Theorem}

\subsection{Hardy spaces and basises}

There are different ways to define Hardy spaces on the the Riemann
surfaces, -- the spaces of vector bundles, multivalued functions or
forms. These definitions are equivalent. We start from
$1$--forms, the most natural object with this respect from our point of
view.

Let $\pi(\zeta):\DD\to X$ be a uniformization of the surface
$X=X_c\setminus E$.  Thus there exists a discrete subgroup $\Gamma$
of the group $SU(1,1)$ consisting of elements of the form
\begin{equation*}
\gamma=\bmatrix \gamma_{11}&\gamma_{12}\\
\gamma_{21}
&\gamma_{22}
\endbmatrix,\
\gamma_{11}=\overline{\gamma_{22}},\
\gamma_{12}=
\overline{\gamma_{21}},
\ \det\gamma=1,
\end{equation*}
such that  $\pi(\zeta)$
is automorphic with respect to $\Gamma$, i.e.,
$\pi(\gamma(\zeta))=\pi(\zeta),\ \forall \gamma
\in\Gamma$,
 and any two preimages of $P\in X$ are
$\Gamma$--equivalent.
  We normalize $Z(\zeta):=(z\circ\pi)(\zeta)$ by
the conditions
$Z(0)=\infty$,
$(\zeta^d Z)(0)>0$.

Note that $\Gamma$
acts dissipatively on
$\TT$ with respect to the Lebesgue measure
$dm$, that is there exists a  measurable
(fundamental) set
$\EE$, which does not
contain any two
$\Gamma$--equivalent points,
and the union
$\cup_{\gamma\in\Gamma}\gamma(\EE)$
is a set of
full measure.
In fact $\EE$ can be chosen as a finite union
of intervals, -- the $\TT$--part of the boundary of the fundamental
domain. For the space of square summable functions  on $\EE$ (with
respect to $dm$), we use
the notation $L^2_{dm|\EE}$.

A character of $\Gamma$ is a complex--valued
function
$\alpha:\Gamma\to \TT$, satisfying
\begin{equation*}
\alpha(\gamma_1\gamma_2)=\alpha(\gamma_1)
\alpha(\gamma_2),\quad\gamma_1,\gamma_2
\in\Gamma.
\end{equation*}
The characters form an Abelian compact
group denoted
by $\Gamma^*$. The further Hardy spaces on $X$ are marked
by characters of $\Gamma$.

Let $f$ be an analytic function in
$\DD$, $\gamma\in\Gamma$.
Then we put
\begin{equation*}
f\vert[\gamma]_k=\frac{f(\gamma(\zeta))}
{(\gamma_{21}\zeta+\gamma_{22})^k}
\quad k=1, 2.
\end{equation*}
Notice that $f\vert[\gamma]_2=
f$ for all $\gamma\in \Gamma$, means
that the form
$f(\zeta)d\zeta$ is invariant with respect to
the substitutions
$\zeta\to\gamma(\zeta)$ ($f(\zeta)d\zeta$ is
an Abelian
integral on $\DD/\Gamma$). Analogously,
$f\vert[\gamma]=
\alpha(\gamma)f$ for all $\gamma\in \Gamma$,
$\alpha\in\Gamma^*$,
means
that the form
$|f(\zeta)|^2\,|d\zeta|$ is invariant with
respect to
these substitutions.

We recall, that a function $f(\zeta)$ is
of Smirnov
class, if it can be represented as a ratio
of two
functions from $H^\infty$ with an outer
denominator. The following spaces related
to the Riemann surface $\DD/\Gamma$ are
counterparts of the standard Hardy spaces $H^2$
($H^1$) on the unit disk.

\begin{definition} \label{def}
The space $A^{2}_1(\Gamma,\alpha)$
($A^{1}_2(\Gamma,\alpha)$) is formed by functions $f$,
which are analytic on $\Bbb D$ and satisfy the
following three conditions
\begin{equation*}
\begin{split}
1)& f \ \text{is of Smirnov class}\\
2)& f\vert[\gamma]=\alpha(\gamma) f\ \ \
(f\vert[\gamma]_2=\alpha(\gamma) f)
\quad
\forall\gamma\in \Gamma\\
3)&
\int_{\Bbb E}\vert f\vert^2\,dm<\infty\  \ \
(\int_{\Bbb E}\vert f\vert\,dm<\infty).
\end{split}
\end{equation*}
\end{definition}

$A^2_1(\Gamma,\alpha)$
is a Hilbert space with the reproducing kernel
$k^\alpha(\zeta,\zeta_0)$, moreover
\begin{equation}\label{widom}
0<\inf_{\alpha\in\Gamma^*}
k^\alpha(\zeta_0,\zeta_0)\le
\sup_{\alpha\in\Gamma^*}
k^\alpha(\zeta_0,\zeta_0)
<\infty.
\end{equation}
Put
\begin{equation*}
k^\alpha(\zeta)=k^\alpha(\zeta,0)\quad\text{and} \quad
K^\alpha(\zeta)=\overline{K^\alpha_\zeta(0)}=\frac{k^\alpha(\zeta)}
{\sqrt{k^\alpha(0)}}.
\end{equation*}

We need one more special function.
  The Blaschke product
\begin{equation*}
b(\zeta)=\zeta
\prod_{\gamma\in\Gamma, \gamma\not= 1_2}
\frac{\gamma(0)-\zeta}{1-\overline{\gamma(0)}
\zeta}
  \frac{\vert\gamma(0)\vert}{\gamma(0)}
\end{equation*}
is called the {\it Green's function}
of $\Gamma$ with
respect to the origin.
It is a character--automorphic function,
i.e., there
exists $\mu\in\Gamma^*$ such that
$b(\gamma(\zeta))=\mu(\gamma)
b(\zeta)$.
  Note, if $G(P)=G(P,P_0)$ denotes
the Green's function of the surface
$X$, then
\begin{equation*}
G(\pi(\zeta))=-\log\vert b(\zeta)\vert.
\end{equation*}

We are ready to construct the basis in $A^2_1(\Gamma,\alpha)$.
Consider the following subspace of this space
$$
\{f\in A^2_1(\Gamma,\alpha):
f(0)=0\}.
$$
The following two facts are evident
\begin{equation*}
\begin{split}
1)&\ \{f\in A^2_1(\Gamma,\alpha):
f(0)=0\}=\{b\tilde f:\tilde f\in A^2_1(\Gamma,\mu^{-1}\alpha)\}=
b A^2_1(\Gamma,\mu^{-1}\alpha),\\
2)&\ A^2_1(\Gamma,\alpha)=
\{K^\alpha\}\oplus
\{f\in A^2_1(\Gamma,\alpha):
f(0)=0\}.
\end{split}
\end{equation*}
Thus
\begin{equation*}
\begin{split}
 A^2_1(\Gamma,\alpha)=&
\{K^\alpha\}\oplus b A^2_1(\Gamma,\mu^{-1}\alpha)\\
=&
\{K^\alpha\}\oplus
\{bK^{\alpha\mu^{-1}}\}\oplus b^2 A^2_1(\Gamma,\mu^{-2}\alpha),
\end{split}
\end{equation*}
and so on.

Basically we proved the following theorem,
note, however, that the second statement is not a direct consequence
of the first one.

\begin{theorem}\label{basis0}
Given $\alpha\in\Gamma^*$, the system of functions $\{b^n
K^{\alpha\mu^{-n}}\}_{n\ge 0}$ forms an orthonormal basis in
$A^2_1(\Gamma,\alpha)$;
the system  $\{b^n
K^{\alpha\mu^{-n}}\}_{n\in\ZZ}$ is an orthonormal basis in
$L^2_{dm|\EE}$.
\end{theorem}

\subsection{The Global Functional Model}
Of course the constructions in this section and our speculation
in Sect. 1 are closely related and of mutual influence. In this
subsection we close the construction by proving the
Global Functional Model Theorem.

Let $\Gamma_0:=\ker\mu$, that is
$\Gamma_0=\{\gamma\in\Gamma:
\mu(\gamma)=1\}$.  Evidently, $b(\zeta)$ and
$(zb^d)(\zeta)$ are holomorphic functions on
the surface
$X_0=\DD/\Gamma_0$.

Assume that $\alpha_0\in\Gamma_0$ can be
extended to a character on $\Gamma$, i.e.,
\begin{equation*}
\Omega_{\alpha_0}=
\{\alpha\in\Gamma^*: \alpha|\Gamma_0=
\alpha_0\}\not=\emptyset.
\end{equation*}
Note that the set of characters
\begin{equation*}
\Omega_{\iota}=
\{\alpha\in\Gamma^*: \alpha|\Gamma_0=
\iota\}
\end{equation*}
where $\iota(\gamma)=1$ for all
$\gamma\in\Gamma_0$ is isomorphic to
the set $(\Gamma/\Gamma_0)^*$.

Let us fix an element
$\hat\alpha_0\in\Omega_{\alpha_0}$. Since
\begin{equation*}
  \{\alpha\in\Gamma^*: \alpha|\Gamma_0=
\alpha_0\} =
\{\hat\alpha_0\beta:\beta\in\Gamma^*:
\beta|\Gamma_0=
\iota\}
\end{equation*}
we can define a measure
$d\chi_{\alpha_0}(\alpha)$ on
$\Omega_{\alpha_0}$ by the relation
\begin{equation*}
d\chi_{\alpha_0}(\alpha)=
d\chi_{\alpha_0}(\hat\alpha_0\beta)
=d\chi_{\iota}(\beta),
\end{equation*}
where $d\chi_{\iota}(\beta)$ is the Haar
measure on $(\Gamma/\Gamma_0)^*$
(the measure $d\chi_{\alpha_0}(\alpha)$
does not depend on a choice of the
element $\hat\alpha_0$).

Obviously,
  $\mathcal T\alpha:=\mu^{-1}\alpha$ is an invertible
ergodic measure--preserving transformation
on $\Omega=\Omega_{\alpha_0}$
with respect to the measure
$d\chi=d\chi_{\alpha_0}$.

The following Theorem is a slightly modified version of
Theorem 2.2 from    \cite{VYu}.

\begin{theorem}\label{model}
 With respect to the basis from Theorem \ref{basis0}, the multiplication
operator by $z$ is  a $2d+1$--diagonal ergodic finite difference
operator with $\Omega=\Omega_{\alpha_0}$,
  $d\chi=d\chi_{\alpha_0}$,
  $\mathcal T\alpha:=\mu^{-1}\alpha$ and
$\alpha_0=\alpha|\Gamma_0$.
Moreover, the operators $\hS_+$ and
$(\hJ\hS^d)_+$ are unitary equivalent to
multiplication by $b$ and $(b^d Z)$
in $A^2_1(\Gamma_0,\alpha_0)$ respectively.
This unitary map is given by the formula
\begin{equation*}
\sum_{\{\gamma\}\in\Gamma/\Gamma_0}
f|[\gamma]\alpha^{-1}(\gamma)=
\sum_{n\in\ZZ_+}x_n(\alpha)b^n
K^{\alpha\mu^{-n}},\ \ f\in A_1^2(\Gamma_0, \alpha_0),
\end{equation*}
where
the vector function
$x(\alpha):=
\{x_n(\alpha)\}$ belongs to $\nLp$.
\end{theorem}

\subsection{Uniqueness Theorem} The natural question up to which extend
our functional realization is unique?
\begin{theorem}
Assume that a finite difference ergodic operator has a finite band
functional model that is there exist a triple $\{X_c, \tilde z, E
 \}$, a character $\alpha_0\in\Gamma^*$ and a map $F$ from $\Omega$
to $\tilde\Omega:=\Omega_{\alpha_0}$ such that
$F\mathcal T\omega=\mu^{-1}F\omega$,  $\chi (F^{-1}(A))=\tilde\chi(A)$,
$A\subset\tilde \Omega$, with $d\tilde\chi:=d\chi_{\alpha_0}$,
here $\mu$ is the character of the Green's function $\tilde b$ on
$X_c\setminus E$. Moreover $q^{(k)}(\omega)=\tilde
q^{(k)}(F\omega)$, where the coefficients $\tilde q^{(k)}(\alpha)$
are generated by the multiplication operator $\tilde z$ with
respect to the orthonormal basis $\{\tilde b^n
K^{\alpha\mu^{-n}}\}_{n\in\ZZ}$.

If the functions $\tilde z$ and $\{d\log\tilde b/d\tilde z\}$
separate points on $X_c\setminus E$ then any local functional
model is generated by one of the branches of the function $\tilde
b$.
\end{theorem}
For proof see \cite{SVYU}.

The following example   shows that in the case when these two
functions $\tilde z$ and $\{d\log\tilde b/d\tilde z\}$ do not
separate points on $X_c\setminus E$ one can give different {\it
global} functional realizations for the same ergodic operator.

\begin{proof}[Example \cite{SVYU}]
Let $J=S^d+S^{-d}$. There exist a "trivial" functional model with
$X_c\setminus E\sim \DD$. In this case $J$ is the multiplication
operator by $z=\zeta^d+\zeta^{-d}$ with respect to the standard
basis $\{\zeta^l\}$ in $L^2_\TT$. Note that $b=\zeta$, thus
$$
w:=\frac{d\log b}{dz}=\frac 1{\zeta^d-\zeta^{-d}}\frac 1 d,
$$
that is $z^2+(w d)^{-2}=4$, $|(wd)^{-1}+z|<2$.

On the other hand let us fix any polynomial $T(u)$, $\deg T=d$,
with real critical values on $\RR\setminus[-2,2]$ and define
$X_c\setminus
E=T^{-1}(\overline\CC\setminus[-2,2])\sim\overline\CC\setminus
T^{-1}[-2,2]$. As we discussed the last set is the resolvent set
for a $d$--periodic Jacobi matrix, say $J_0$. Moreover
$T(J_0)=J$, and $-\log|b|$ is just the Green's function of this
domain in the complex plain. So, using the standard functional
model for $J_0$ (see Sect. 7 for details) with the symbols $u$ and $b$
we get a functional model for $J$ with $z=T(u)$ and the same $b$. Note
that as before
$z^2+(w d)^{-2}=4$, $|(wd)^{-1}+z|<2$ with $w:=\frac{d\log
b}{dz}$.
\end{proof}

\begin{remark} In \cite{DKS} the identity $T(J_0)=S^d+S^{-d}$ that
holds for a periodic Jaconi  matrix $J_0$ with the spectrum on
$T^{-1}[-2,2]$ is called the Magic Formula. There it plays an
important role in proving counterparts of Denisov--Rakhmanov and
Killip--Simon Theorems for perturbations of periodic Jacbi matrices.
\end{remark}


\section{Several infinities case. Existence Theorem}

Now we examine the situation
in which the reduced surface $X_c$
has several "infinities" that is the covering function $z$
(the symbol of an almost periodic operator) equals infinity
at several (distinct) points on $X_c$.

We start with a simple example.

\subsection{A five diagonal matrix of period two}
Assume that $z$ is a two sheeted covering with only two branching
points, say, $z_1=-2$, $z_2=2$.
 With necessity corresponding substitutions are $\sigma_1=\sigma_2=
\left(\begin{matrix} 1& 2\\2&1\end{matrix}\right)$. In this case
$X_c$ is equivalent to the complex plane
$\overline\CC$, moreover we can put
\begin{equation}\label{vs1}
z= v+\frac 1 v,\quad v\in \CC.
\end{equation}
Thus, on this surface we have two "infinities" $v=\infty$ and
$v=0$.

Note that $z^{-1}(\RR)=\TT\cup\RR$.
We cut $\overline \CC$ over the interval $[a,b]$, where
$0<a<b<1$, that is we consider $X_c\setminus E$ of the form
$\overline \CC\setminus [a,b]$. We uniformize
$X_c\setminus E$ by
\begin{equation}\label{vs2}
v=\frac{a+b}{2}+\frac{b-a}{2}\frac{1/\zeta+\zeta}{2},
\quad \zeta\in \DD.
\end{equation}
 For $v=\infty$ we have $\zeta_0=0$.
Solving
$$
\zeta^2+2\frac{b+a}{b-a}\zeta+1=0
$$ we get the image of the second infinity in $\DD$,
$\zeta_1=-\frac{\sqrt{b}-\sqrt{a}}{\sqrt{b}+\sqrt{a}}$.
As a result we get a symbol function $z:\DD\to \CC$
(see \eqref{vs1}, \eqref{vs2}) for the forthcoming operator
with two infinities $\zeta_0, \zeta_1$.

Next point is the symbol $b$ for the shift operator. Recall that the
Green function in
$\DD$ is the Blaschke factor
$$
b_{\zeta_0}=\zeta,\quad
b_{\zeta_1}=\frac{\zeta-\zeta_1}{1-\zeta\zeta_1}.
$$
The product $b_{\zeta_0}b_{\zeta_1}$ is the smallest unimodular
multiplier that cancels poles of $z$.
Since $b^2$ is of the same nature, $b^2 z$ is holomorphic with a
unimodular function
$b$ on $\TT$,
we have $b^2=b_{\zeta_0}b_{\zeta_1}$.

Finally, we need a certain functional space
and an intrinsic basis in it
that  generalize the construction in Theorem \ref{basis0}.
Recall $b$ should be related
to the shift $S$, and we are going to define the periodic
operator
$J$ as the multiplication operator with respect to this basis.
To this end we define
the following functional spaces. Given $\alpha_k\in \TT$,
$k=0,1$, we associate the space
$H^2(\alpha_0,\alpha_1)$ of
analytic multivalued functions $f(\zeta)$,
$\zeta\in\DD\setminus\{\zeta_0,\zeta_1\}$, such that
$|f(\zeta)|^2$ has a harmonic majorant and
$$
f\circ\gamma_i=\alpha_i f,
$$
where $\gamma_i$ is a small circle around $\zeta_i$.
Such a space can be reduced to the standard Hardy space $H^2$, moreover
$$
H^2(\alpha_0,\alpha_1)=
b_{\zeta_0}^{\tau_0}b_{\zeta_1}^{\tau_1}H^2,
\quad \alpha_k=e^{2\pi i\tau_k},\ 0\le\tau_k<1.
$$

\begin{lemma} The space $bH^2(-1,1)$ is a subspace
of $H^2(1,-1)$ having a one  dimensional orthogonal compliment,
moreover
\begin{equation}\label{vs3}
H^2(1,-1)=\{\sqrt{b_{\zeta_1}}k_{\zeta_0}\}
\oplus bH^2(-1,1),
\end{equation}
where $k_{\zeta_i}$ is the reproducing kernel of the standard $H^2$
with respect to $\zeta_i$.
\end{lemma}

This lemma allows us to repeat construction of subsection 3.1.
Iterating, now, the decomposition \eqref{vs3}
\begin{equation*}
\begin{split}
H^2(1,-1)=&\{\sqrt{b_{\zeta_1}}k_{\zeta_0}\}
\oplus bH^2(-1,1)\\=&
\{\sqrt{b_{\zeta_1}}k_{\zeta_0}\}
\oplus
b\{\sqrt{b_{\zeta_0}}k_{\zeta_1}\}
\oplus
b^2H^2(1,-1)=...,
\end{split}
\end{equation*}
one gets an orthogonal basis
in $H^2(1,-1)$ consisting of vectors of two sorts
$$
b^{2m}\{\sqrt{b_{\zeta_1}}k_{\zeta_0}\}\ \ \text{and}\ \
b^{2m+1}\{\sqrt{b_{\zeta_0}}k_{\zeta_1}\}.
$$
Note that this orthogonal
system can be extended on negative integers $m$ so that we obtain
a basis in $L^2(1,-1)$.

\begin{theorem} With respect to the orthonormal basis
\begin{equation}\label{ts4}
e_n=\begin{cases}
b^{2m}\sqrt{b_{\zeta_1}}\frac{k_{\zeta_0}}{||k_{\zeta_0}||}
, &n=2m\\
b^{2m+1}\sqrt{b_{\zeta_0}}\frac{k_{\zeta_1}}{||k_{\zeta_1}||}
,&n=2m+1
\end{cases}
\end{equation}
the multiplication operator by $z$ is a 5--diagonal  matrix of
period 2.
\end{theorem}

\subsection{General case}
Let $z: X_c\to \CC$ be $d$--sheeted covering with $d$ (distinct)
infinities. Further, let $\DD/\Gamma$ be a uniformization
of $X_c\setminus E$. A given character $\alpha\in \Gamma^*$ and a
{\it fixed ordering}
of infinities $P_1, P_2,...,P_d$ define
\begin{itemize}
\item $b_j=b_{P_j}$ the Green function
with respect to $P_j$, $b_j\circ\gamma=\mu_j(\gamma) b_j$,
$\mu_j\in\Gamma^*$;
\item $k_j^\alpha=k_{P_j}^\alpha$ the reproducing kernel of
$A_1^2(\Gamma,\alpha)$ with respect to $P_j$,
$K_j^\alpha:=\frac{k_j^\alpha}{||k_j^\alpha||}$.
\end{itemize}
\begin{theorem} \label{t1.3}
Let $b=(b_1...b_d)^{\frac 1 d}$.
With respect to the orthonormal basis
\begin{equation}\label{ts5}
e_n=\begin{cases}
b^{dm}b_1^{\frac{1}{d}}...b_{d-1}^{\frac{d-1}{d}}
K_d^{\alpha(\mu_1...\mu_d)^{-m}}
, &n=dm\\
b^{dm+1}b_2^{\frac{1}{d}}...b_{d}^{\frac{d-1}{d}}
K_1^{\alpha\mu_1^{-1}(\mu_1...\mu_d)^{-m}}
, &n=dm+1\\
\dots&\\
b^{dm+d-1}b_d^{\frac{1}{d}}...b_{d-2}^{\frac{d-1}{d}}
K_{d-1}^{\alpha\mu_1^{-1}...\mu_{d-1}^{-1}(\mu_1...\mu_d)^{-m}}
, &n=dm+d-1
\end{cases}
\end{equation}
the multiplication operator by $z$ is a $(2d+1)$--diagonal
almost periodic matrix.
\end{theorem}

\section{Five--diagonal almost periodic self--adjoint matrices
and OPUC}

We start again with two--sheeted covering \eqref{vs1}.
We have $X_c\simeq \overline\CC$ and $z^{-1}(\RR)=\RR\cup\TT$,
but let us, in this case, cut $\overline\CC$
on an arbitrary finite union of (necessary non--degenerate) arcs
on the unit circle,
$E\simeq\{\TT\setminus\cup_{j=0}^l(a_j,b_j)\}$.

The domain $X_c\setminus E$
is conformally equivalent to the quotient of the
unit disk by the action of a discrete group $\Gamma=\Gamma(E)$.
Let
\begin{equation}\label{8s6}
v:\DD\to \{\overline\CC\setminus\TT\}\cup\{\cup_{j=0}^l(a_j,b_j)\}
\end{equation}
be a covering map,
$v\circ\gamma=v$, $\forall \gamma\in\Gamma$. In what follows we
assume the following normalization to be hold
$
v:(-1,1)\to (a_0,b_0),
$
so that one can chose
a fundamental domain $\mathfrak F$ and a system of generators
$\{\gamma_j\}_{j=1}^l$ of $\Gamma$ such that
they are symmetric with respect to the complex conjugation:
$$
\overline{\mathfrak F}=\mathfrak F, \quad
\overline{\gamma_j}=\gamma_j^{-1}.
$$
Denote by $\zeta_0\in\mathfrak F$ the preimage of the origin,
$v(\zeta_0)=0$, then $v(\overline{\zeta_0})=\infty$.
Thus infinities of $z$ defined by \eqref{vs1} and \eqref{8s6}
on $\DD/\Gamma$ are trajectories of $\zeta_0$ and $\overline{\zeta_0}$
under the action of the group $\Gamma$,
$P_0=\{\gamma(\zeta_0)\}_{\gamma\in \Gamma}$
and $P_1=\{\gamma(\overline{\zeta_0})\}_{\gamma\in \Gamma}$.
These are {\it two} infinities that we have in the case under
consideration.

Thus, to define the function $b$ (the symbol of the shift operator)
we have to introduce the Green functions
$B(\zeta,\zeta_0)$ and
$B(\zeta,\overline{\zeta_0})$. It is convenient to
normalize them by
$B(\overline{\zeta_0},\zeta_0)>0$ and
$B(\zeta_0,\overline{\zeta_0})>0$. Then
\begin{equation}\label{8s7}
v(\zeta)=e^{ic}\frac{B(\zeta,{\zeta_0})}{B(\zeta,\overline{\zeta_0})}.
\end{equation}
Also, we can rotate (if necessary) $v$--plane so
that $c=0$. Note that $B(\zeta,{\zeta_0})$ is a
character--automorphic function
$$
B(\gamma(\zeta),{\zeta_0})=\mu(\gamma)
B(\zeta,{\zeta_0}), \quad\gamma\in\Gamma,
$$
with a certain $\mu\in\Gamma^*$. By \eqref{8s7},
$B(\zeta,\overline{\zeta_0})$ has the {\it same} factor of automorphy,
$$
B(\gamma(\zeta),\overline{\zeta_0})=\mu(\gamma)
B(\zeta,\overline{\zeta_0}), \quad\gamma\in\Gamma.
$$
By the definition $b^2=B(\zeta,{\zeta_0})
B(\zeta,\overline{\zeta_0})$,
we get a multivalued analytic function $b$
on the punched surface $\{\DD/\Gamma\}\setminus\{P_0,P_1\}$.

In this case we have only two possibilities for ordering of infinities:
$\{P_0,P_1\}$ and $\{P_1,P_0\}$. According to Theorem \ref{t1.3},
to any of them, say the first one,  and to an arbitrary
$\alpha\in{\Gamma^*}$ we can associate
the operator $J=J(-1,1;\alpha)$ by fixing the space
$H^2(-1,1;\alpha)=\sqrt{B(\zeta,\zeta_0)}A_1^2(\alpha)$ and
a natural basis in it. Up to a common multiplier
$\sqrt{B(\zeta,\zeta_0)}$
it is a basis in
$A_1^2(\alpha)$ of the form
$$
K^{\alpha}(\zeta,\overline{\zeta_0}),\,
B(\zeta,\overline{\zeta_0})K^{\alpha\mu^{-1}}(\zeta,{\zeta_0}),\,
B(\zeta,{\zeta_0})
B(\zeta,\overline{\zeta_0})K^{\alpha\mu^{-2}}(\zeta,\overline{\zeta_0}),
\,...
$$

Thus we get the same system of functions that we used
describing almost periodic Verblunsky coefficients
\cite{PYU}. The last one we can define by
$$
a(\alpha)=a=\frac{K^{\alpha}(\zeta_0,\overline{\zeta_0})}
{K^{\alpha}(\zeta_0,{\zeta_0})}.
$$
In \cite{PYU} they appear in
 the following recursion
\begin{equation}\label{8s8}
\begin{matrix}
K^{\alpha}(\zeta,\overline{\zeta_0})=&
a(\alpha)K^{\alpha}(\zeta,{\zeta_0})
+\rho(\alpha)
B(\zeta,\zeta_0)K^{\alpha\mu^{-1}}(\zeta,\overline{\zeta_0}),
\\
K^{\alpha}(\zeta,{\zeta_0})=&\overline{a(\alpha)}
K^{\alpha}(\zeta,\overline{\zeta_0})
+\rho(\alpha)
B(\zeta,\overline{\zeta_0})K^{\alpha\mu^{-1}}(\zeta,{\zeta_0}),
\end{matrix}
\end{equation}
where
$$
\rho(\alpha)=\rho=\sqrt{1-|a|^2}=
B(\overline{\zeta_0},\zeta_0)
\frac{K^{\alpha\mu^{-1}}(\overline{\zeta_0},\overline{\zeta_0})}
{K^{\alpha}(\overline{\zeta_0},\overline{\zeta_0})}.
$$

The goal of this section is to represent $J$ in terms of
Verblunsky coefficients.
\begin{lemma}
With respect to the basis
\begin{equation}\label{bas}
...,\,
K^{\alpha\mu}(\zeta,{\zeta_0})/B(\zeta,{\zeta_0}),\,
K^{\alpha}(\zeta,\overline{\zeta_0}),\,
B(\zeta,\overline{\zeta_0})K^{\alpha\mu^{-1}}(\zeta,{\zeta_0}),
\,...
\end{equation}
the multiplication operator by $v$ is a matrix having at most
two non vanishing diagonals over the main diagonal. Moreover,
\begin{equation}\label{9s9}
v\sim
\begin{bmatrix}
 &\rho(\alpha\mu) \rho(\alpha\mu^2)  & 0 & \\
\ddots &-\rho(\alpha\mu)  a(\alpha\mu^2) & 0 & \\
&  -a(\alpha\mu)  \overline{a(\alpha)}& \rho(\alpha)
\overline{a(\alpha\mu^{-1})}&
\\ &  *& -a(\alpha) \overline{a(\alpha\mu^{-1})}& \\
& & & \ddots\\
\end{bmatrix}.
\end{equation}

Similarly, the multiplication operator by $1/v$ is of the form
\begin{equation}\label{9s10}
1/v\sim
\begin{bmatrix}
&0  &0 & \\
\ddots &\rho(\alpha\mu)  a(\alpha) & \rho(\alpha)
\rho(\alpha\mu)  & \\
&  -\overline{a(\alpha\mu)} a(\alpha)&
-\rho(\alpha) \overline{a(\alpha\mu)}&
\\ &  *& -a(\alpha\mu^{-1}) \overline{a(\alpha)}& \\
& & & \ddots\\
\end{bmatrix}.
\end{equation}
\end{lemma}

\begin{proof} We give a proof, say, for \eqref{9s9}.
Recall \eqref{8s7}, from which we can see that
the decomposition of the vector
$v(\zeta) K^{\alpha}(\zeta,\overline{\zeta_0})$ begins with
$$
v(\zeta) K^{\alpha}(\zeta,\overline{\zeta_0})=
c_0\frac{K^{\alpha\mu^2}(\zeta,\overline{\zeta_0})}
{B(\zeta,{\zeta_0})B(\zeta,\overline{\zeta_0})}+
c_1\frac{K^{\alpha\mu}(\zeta,{\zeta_0})}
{B(\zeta,{\zeta_0})}+
c_2{K^{\alpha}(\zeta,\overline{\zeta_0})}
+ ...\,.
$$
 Multiplying by the denominator
$B(\zeta,{\zeta_0})B(\zeta,\overline{\zeta_0})$ we get
\begin{equation}\label{9s11}
\begin{split}
B^2(\zeta,{\zeta_0}) K^{\alpha}(\zeta,\overline{\zeta_0})=&
c_0{K^{\alpha\mu^2}(\zeta,\overline{\zeta_0})}
+
c_1{K^{\alpha\mu}(\zeta,{\zeta_0})}
{B(\zeta,\overline{\zeta_0})}\\+&
c_2{K^{\alpha}(\zeta,\overline{\zeta_0})}
{B(\zeta,{\zeta_0})B(\zeta,\overline{\zeta_0})}
+ ...\,.
\end{split}
\end{equation}
First we put $\zeta=\overline{\zeta_0}$.
By the definition of $\rho(\alpha)$ we have
$$
c_0=B^2(\overline{\zeta_0},{\zeta_0})
\frac{K^{\alpha}(\overline{\zeta_0},\overline{\zeta_0})}
{K^{\alpha\mu^2}(\overline{\zeta_0},\overline{\zeta_0})}
=\rho(\alpha\mu) \rho(\alpha\mu^2).
$$
Putting $\zeta={\zeta_0}$ in \eqref{9s11}
and using the definition of $a(\alpha)$,
we have
$$
c_1=-c_0\frac{K^{\alpha\mu^2}(\zeta_0,\overline{\zeta_0})}
{K^{\alpha\mu}(\zeta_0,{\zeta_0})
B(\zeta_0,\overline{\zeta_0})}=
-\rho(\alpha\mu) \rho(\alpha\mu^2)
\frac{a(\alpha\mu^2)}{\rho(\alpha\mu^2)}=
-\rho(\alpha\mu) a(\alpha\mu^2).
$$

Doing in the same way we can find a representation for $c_2$
that would involve derivatives of the reproducing kernels. However,
we can find
$c_2$ in terms of
$a$ and
$\rho$  calculating the scalar product
$$
c_2=\langle
B^2(\zeta,{\zeta_0}) K^{\alpha}(\zeta,\overline{\zeta_0}),
{B(\zeta,{\zeta_0})B(\zeta,\overline{\zeta_0})}
{K^{\alpha}(\zeta,\overline{\zeta_0})}
\rangle.
$$
Since $B(\zeta,{\zeta_0})$ is unimodular, using
\eqref{8s8}, we get
$$
c_2=\langle\frac{K^{\alpha\mu}(\zeta,\overline{\zeta_0})
- a(\alpha\mu)K^{\alpha\mu}(\zeta,\zeta_0)}
{\rho(\alpha\mu)},
{B(\zeta,\overline{\zeta_0})}
{K^{\alpha}(\zeta,\overline{\zeta_0})}
\rangle.
$$
Recall that
$
k^{\alpha}(\zeta,\zeta_0)=
K^{\alpha}(\zeta,\zeta_0)K^{\alpha}(\zeta_0,\zeta_0)
$ is the reproducing kernel. Thus
$$
c_2=-\frac{a(\alpha\mu)}{\rho(\alpha\mu)}\overline{
\frac{{B(\zeta_0,\overline{\zeta_0})}
{K^{\alpha}(\zeta_0,\overline{\zeta_0})}}
{K^{\alpha\mu}(\zeta_0,\zeta_0)}}=
-\frac{a(\alpha\mu)}{\rho(\alpha\mu)}
\overline{\rho(\alpha\mu)a(\alpha)}
=-a(\alpha\mu)\overline{a(\alpha)}.
$$

To find the decomposition of the vector
$v(\zeta)B(\zeta,\overline{\zeta_0})
{K^{\alpha\mu^{-1}}(\zeta,{\zeta_0})}$
is even simpler
since only two leading terms are involved.
Note that all other columns of the matrix in
\eqref{9s9}, starting from these two, can be obtain by
the character's shift by $\mu^{-2}$ along diagonals.

\end{proof}

Now let us remind the CMV representation for operators related
to OPUC \cite{S}. A given sequence of numbers from $\DD$
\begin{equation}\label{16s13}
...,\,a_{-1},\,a_{0},\,a_{1},\,a_{2},\,...
\end{equation}
define unitary matrices
\begin{equation*}
A_k=\begin{bmatrix} \overline{a_k}& \rho_k\\
\rho_k&-a_k
\end{bmatrix}, \quad \rho_k=\sqrt{1-|a_k|^2},
\end{equation*}
and unitary operators in $l^2(\ZZ)$ given by
block--diagonal matrices
\begin{equation*}
\mathfrak A_0=\begin{bmatrix} \ddots& & &\\
& A_{-2}& & \\
& & A_{0}&  \\
 & & &\ddots
\end{bmatrix}, \quad
\mathfrak A_1=S\begin{bmatrix} \ddots& & &\\
& A_{-1}& & \\
& & A_{1}&  \\
 & & &\ddots
\end{bmatrix}S^{-1}.
\end{equation*}
The CMV operator $\mathfrak A$, related to
the sequence \eqref{16s13}, is the
product
\begin{equation}\label{16s14}
\mathfrak A=\mathfrak A(\{a_k\}):=\mathfrak A_0\mathfrak A_1.
\end{equation}

\begin{theorem} Define the sequence
$a_k=a(\alpha\mu^{-k})$. Then
$v\sim \mathfrak A(\{a_k\})$, see \eqref{16s14}.
\end{theorem}

\begin{proof}
Note that operators in \eqref{9s9} and \eqref{9s10}
are mutually conjugated, therefore under--diagonal
entries of both operators are also known. The rest is an easy direct
computation.
\end{proof}

\begin{theorem} With respect to the basis \eqref{bas}
the multiplication operator by $z$,
\begin{equation*}
z\sim
\begin{bmatrix}
 &q^{(-2)}(\alpha)   & 0 & \\
\ddots &q^{(-1)}(\alpha)  & q^{(-2)}(\alpha\mu^{-1}) & \\
&  q^{(0)}(\alpha) & \{\overline{q^{(-1)}(\alpha\mu^{-1})}\}  &\\
&  *& q^{(0)}(\alpha\mu^{-1}) & \\ & & & \ddots\\
\end{bmatrix},
\end{equation*}
is defined
by functions
$$
q^{(-2)}(\alpha)=\rho(\alpha\mu) \rho(\alpha\mu^2),\quad
q^{(-1)}(\alpha)=\rho(\alpha\mu)\{a(\alpha)- a(\alpha\mu^2)\},
$$
and
$$
q^{(0)}(\alpha)=-2\Re\{a(\alpha) \overline{a(\alpha\mu)}\}.
$$
\end{theorem}
\noindent
It is worth to mention that the second column is not only the shift
by $\mu^{-1}$, there is also the conjugation.

\begin{proof} Recall that $z=v+1/v$ and use the previous lemma.
\end{proof}

Using general constructions from \cite{MM} one can define
and integrate the flows hierarchy given by
\begin{equation}\label{16s15}
\dot \fA=[(\fA^n+\fA^{-n})_+,\fA].
\end{equation}
Note that for $n=1$, \eqref{16s15} gives the Schur flow
\cite{FG}.

\section{Coverings }

Discussing this subject we prefer to use a functional version of
presentation of Hardy spaces.
First
we introduce these spaces and then will remark how they are related
to the spaces of forms.

Let $\omega_{P_0}$ denote the harmonic measure
on an open surface $X_c\setminus E$ with respect to
$P_0\in X_c\setminus E$. Note that $\omega_{P_0}$ is the restriction
of the differential $\frac 1{2\pi i}d\log b(P,P_0)$ on $E$.
By $H^2(\alpha,\omega_{P_0})$ we denote the closure
of $H^\infty(\alpha)$ in $L^2$ with respect to the measure
$\omega_{P_0}$. The natural question is how this space
is related to the space with another point fixed, say
$P_1\in X_c\setminus E$,
or, more generally, with
$H^2(\alpha,\omega)$, which denotes the closure of $H^\infty(\alpha)$
with respect to an equivalent norm given by a measure of the form
$\omega=\rho
\omega_{P_0}$,  where $0<C_1\le\rho\le C_2<\infty$. (By the Harnack
Theorem $\omega_{P_1}$ and $\omega_{P_0}$ satisfy this property).

To answer it, let us define
an outer function $\phi$,  such that $\rho=|\phi|^2$.
This function belongs to $H^\infty(\beta)$ with a certain
$\beta\in\Gamma^*$
In this case
$$
f\mapsto\phi f
$$
is a unitary map from $H^2(\alpha,\omega)$ to
$H^2(\alpha\beta,\omega_{P_0})$. Then the equality
\begin{equation*}
\begin{split}
\langle (\phi f)(P),
k^{\alpha}_Q(P;\omega)\phi(P)\overline{\phi(Q)}
\rangle_{H^2(\alpha\beta,\omega_{P_0})}=&
\langle f(P),
k^\alpha_Q(P;\omega)\overline{\phi(Q)}
\rangle_{H^2(\alpha,\omega)}\\
=&f(Q)\phi(Q)
\end{split}
\end{equation*}
shows that the reproducing kernels of  $H^2(\alpha\beta,\omega_{P_0})$
and  $H^2(\alpha,\omega)$ are related by
$$
k^{\alpha\beta}_Q(P;\omega_{P_0})=
k^{\alpha}_Q(P;\omega)\phi(P)\overline{\phi(Q)}.
$$
For the normalized kernels we have
$$
K^{\alpha\beta}_Q(P;\omega_{P_0})=
K^{\alpha}_Q(P;\omega)\frac{\overline{\phi(Q)}}{|\phi(Q)|}
\phi(P).
$$
Therefore the matrix of a multiplication operator
in $H^2(\alpha,\omega)$,
with respect to the reproducing kernels  based basis,
 actually  can be obtain by a character's shift
for the matrix with the same symbol related to the
 chosen space $H^2(\alpha,\omega_{P_0})$.
 (Let us mention in brackets a specific normalization
 of a basis vector
 given by the unimodular factor
$\frac{\overline{\phi(Q)}}{|\phi(Q)|}$).

The relations between $A_1^2(\alpha)$ and $H^2(\alpha)$
are of the same nature. Indeed, let
$\rho:\DD/\Gamma\to X_c\setminus E$, $\rho(0)=P_0$, be the
uniformization of the given surface. Then $H^2(\alpha)$ is a subspace of
the standard $H^2$ in $\DD$: $f(\rho(\zeta))\in H^2$
for $f\in H^2(\alpha)$, moreover
$$
\Vert f\Vert^2=\int_{\TT}|f(\rho(t))|^2\,dm(t),
$$
where $dm$ is the Lebesgue measure. Fix a fundamental set $\EE$
for the action of $\Gamma$ on $\TT$,
$\TT=\cup_{\gamma\in\Gamma}\gamma(\EE)$. Then
$$
\Vert f\Vert^2=\int_{\EE}|f(\rho(t))|^2|\psi(t)|^2\,dm(t),
$$
where
$$
|\psi(t)|^2:=\sum_{\gamma\in\Gamma}|\gamma'(t)|.
$$
Again, we can consider $\psi$ as an outer function and then
$$
f\to  (f\circ\rho)  \psi
$$
is the unitary map from  $H^2(\alpha)$ to
$A_1^2(\alpha\beta)$, where the character $\beta$
is generated by the 1--form $\psi$.
The $\beta$ here is a particular character, so when
$\alpha$ runs on the whole group $\Gamma^*$
$\alpha\beta$ covers also all characters, and we have one to one
correspondence between two ways of writing of the Hardy spaces.

But, as we noted above, working with coverings, it will be convenient to
use character automprphic
$H^2$--spaces with respect to the following specific measure
\begin{equation}\label{24s16}
\omega=\frac 1 l \sum_{i=1}^l \omega_{P_i},
\end{equation}
associated with a system of points $\{P_i\}$ on $X_c\setminus E$.
Naturally in what follows $P_i$'s are infinities on $X_c\setminus E$.

Now we can go back to the coverings.
Let $\Gamma_X$ (respectively $\Gamma_Y$) be the fundamental
group on $X_c\setminus E$ (respectively $Y_c\setminus E$).
We have $\pi_* :\Gamma_Y\to \Gamma_X$ ($\pi_*(\gamma)$ is the image
of a contour $\gamma\in \Gamma_Y$) and
$\pi^* :\Gamma_X\to \Gamma_Y$ ($\pi^*(\gamma)$ is the full preimage
of a contour $\gamma\in \Gamma_X$). Note that
$\pi^*\pi_*=d\, {\rm Id}$. The maps $\pi_* :\Gamma_Y^*\to \Gamma_X^*$
and $\pi^* :\Gamma_X^*\to \Gamma_Y^*$ are defined by duality.

For a system of points $\{P_i\}_{i=1}^l$, $P_i\in X_c\setminus E$,
define the measure $\omega_*$ on $E$  by \eqref{24s16}.
Let $\{Q^{(i)}_k\}_{k=1}^d=\pi^{-1}(P_i)$. Define
\begin{equation*}
\omega^*=\frac 1 {ld} \sum \omega_{Q^{(i)}_k}.
\end{equation*}
In this case
\begin{equation}\label{24s17}
\int_F (f\circ\pi)\omega^*=\int_E f\omega_*.
\end{equation}
Moreover,
\begin{equation}\label{24s18}
\int_F f\omega^*=\int_E (\fL f)\omega_*,
\quad
{\rm where}\quad
(\fL f)(P)=\frac 1 d \sum_{\pi(Q)=P}f(Q).
\end{equation}

As a direct consequence of \eqref{24s17}, \eqref{24s18}
we get
\begin{lemma} The map $V: H^2(\alpha,\omega_*)\to
H^2(\pi^*\alpha,\omega^*)$, defined by
\begin{equation}\label{24s19}
V f=f\circ\pi,\quad f\in H^2(\alpha,\omega_*),
\end{equation}
is an isometry with
\begin{equation}\label{24s20}
(V^* f)(P)=(\fL f)(P), \quad f\in
H^2(\pi^*\alpha,\omega^*).
\end{equation}
Also,
\begin{equation}\label{24s20}
V k^\alpha_{P_0}=\frac 1 d \sum_{\pi(Q_j^{(0)})=P_0}
k^{\pi^*\alpha}_{Q_j^{(0)}}.
\end{equation}
for the reproducing kernel  $k^\alpha_{P_0}\in H^2(\alpha,\omega_*)$.
\end{lemma}

\begin{theorem} Let $z_*:X_c\setminus E\to \overline{\CC}$.
Using notations introduced above, assume that
$$
z_*: E\to \RR, \quad z_*^{-1}(\infty)\subset{\{P_i\}_{i=1}^l}
$$
Let $\mathbf z_*(\alpha)$ be the multiplication operator
by $z_*$ with respect to the basis \eqref{ts5}.
For an arbitrary ordering of $\{Q^{(i)}_k\}_{1\le k\le d, 1\le i\le l}$
subordinated to the ordering of $\{P_i\}$ consider the operator
$\mathbf z^*(\pi^*\alpha)$ with the symbol $z^*:=z_*\circ\pi$
and the related by \eqref{24s19} isometry $\mathbf V$.
Then, the following relations hold
\begin{equation}\label{24s}
\mathbf V^*\mathbf z^*(\pi^*\alpha)=\mathbf z_*(\alpha)\mathbf V^*,
\quad
\mathbf V^*S^d=S\mathbf V^*.
\end{equation}
\end{theorem}

Equations \eqref{24s} are very close to the so called
 Renormalization Equations that we start to discuss now.


\section{
The renormalization
of periodic matrices}

We recall some basic facts from the spectral theory of periodic
Jacobi matrices. The spectrum $E$ of any periodic matrix $J$ is an
inverse polynomial image
\begin{equation}\label{4.4}
E=U^{-1}[-1,1]
\end{equation}
the polynomial $U$ of degree $g+1$ should have
all critical points $\{c_U\}$ real and for all critical values
$|U(c_U)|\ge 1$. For simplicity we assume
that
$|U(c_U)|> 1$.
Then the spectrum of $J$ consists of $g$ intervals
$$
E=[b_0,a_0]\setminus (\cup_{j=1}^g(a_j,b_j)).
$$
Also it would be convenient for us to normalize $U$ by a linear
change of the variable such that $b_0=-1$ and $a_0=1$.
In this case $U$ is a so called expanding polynomial.

Having the set $E$ of the above form fixed, let us describe the whole set
of periodic Jacobi matrices $J(E)$ with the given spectrum.
To this end we associate with $U$ the hyper--elliptic Riemann surface
(the surface is given by \eqref{5.06} with $\lambda =b^N$,
$N=g+1$)
$$
X=\{Z=(z,\lambda): \lambda-2 U(z)+\lambda^{-1}=0\}.
$$
The involution on it we denote by $\tau$,
\begin{equation}\label{invol}
\tau Z:=\left(z,\frac{1}{\lambda}\right)\in X.
\end{equation}
The set
$$
X_+=\{Z\in X: |\lambda(Z)|<1\}
$$
we call the upper sheet of $X$. Note $X_+\simeq \bar \CC\setminus E$,
in fact, $z(Z)\in \bar \CC\setminus E$ if $Z\in X_+$.

The following well known theorem describes $J(E)$ in terms of {\it real}
divisors on
$X$. The Jacobian variety of $X$, $\Jac(X)$, is a $g$ dimensional
complex torus,
$\Jac(X)\simeq \CC^g/L(X)$, where $L$ is a lattice (that can be chosen
in the form $L=\ZZ^g+\Omega \ZZ^g$ with $\Im \Omega>0$). Consider the
$g$ dimensional
real subtorus consisting of divisors of the form
$$
D(E)=\{D=D_+ -D_C,\ D_+:= \sum_{i=1}^g Z_i: Z_i\in X,
\ z(Z_i)\in[a_i,b_i]\},
$$
here $D_C$ is a point of normalization that we choose of the form
$$
D_C:=\sum_{i=1}^g C_i
: C_i\in X,
\ z(C_i)=(c_U)_i, |\lambda(C_i)|>1,
$$
--- the collections of the points on the lower sheet with the
$z$--coordinates at the critical points. (At least topologically,
it is evident $D(E)\simeq \RR^g/\ZZ^g$).

\begin{theorem}\label{param}
For given $E$ of the form \eqref{4.4} there exists an one--to--one
correspondence
between $J(E)$ and $D(E)$.
\end{theorem}

Let now $\tilde U$ and $T$ be polynomials of the described  above form,
and we define $U=\tilde U\circ T$. Then we have a covering
$\pi$ of the Riemann surface
$\tilde X$ associated to $\tilde U$ by the surface $X$ associated to $U$:
\begin{equation}\label{cov}
\pi(z,\lambda)=(T(z), \lambda),
\end{equation}
note $\pi: X_+\to \tilde X_+$.

According to the general theory, this covering generates different natural
mappings \cite{Mam},
in particular,
\begin{equation}\label{covjac}
\pi_*: \Jac(X)\to \Jac(\tilde X),
\end{equation}
and
\begin{equation}\label{covjacback}
\pi^*: \Jac(\tilde X)\to \Jac(X).
\end{equation}
Thus, in combination with Theorem \ref{param},
we get the map
\begin{equation}\label{themap}
\begin{array}{lll}
D(\tilde E)
&\stackrel{\pi_T^*}{\longrightarrow}& D(E)\\
\ \ \Big\downarrow & &
\ \ \Big\downarrow\\
J(\tilde E)&\stackrel{{\mathcal J}_T}{\longrightarrow}& J(E)
\end{array}
\end{equation}
We study this map in terms  described in the previous section.

\bigskip


The differential $ \frac{1}{2\pi i}d\log b$,
being restricted
on $\partial X_+$, is the  harmonic measure $d\omega$ of the domain
$\bar
\CC\setminus E$
with pole at infinity.
The  space $L^p(\partial X_+)$, in a sense, is the $L^p$
space with respect to the harmonic measure, but it should be mentioned
that $\partial X_+= (E-i0)\cup(E+i0)$, i.e.,
an element $f$ of $L^p(\partial X_+)$ may have different values
$f(x+i0)$ and $f(x-i0)$, $x\in E$,
and $H^p(X_+)$ is the closure of the set of
holomorphic functions uniformly bounded in $ X_+$ with respect
to this norm.

Since $\tilde \lambda\circ \pi= \lambda$ we have the relation
\begin{equation}\label{4.7s}
\tilde b\circ\pi= b^d.
\end{equation}
Thus, we get
\begin{equation}\label{4.8s}
\int_{\partial X_+} f\, d\omega
= \int_{\partial \tilde X_+} \frac 1 d
\left(\sum_{\pi( Z)= \tilde Z}f(Z) \right)(\tilde Z)\, d\tilde\omega
\end{equation}
for every $f \in L^1(\partial X_+)$.

As it follows directly from \eqref{4.8s},
the covering \eqref{cov} generates an isometrical enclosure
\begin{equation}\label{encl1}
v_+: H^2(\tilde X) \to H^2(X_+)
\end{equation}
acting in a natural way
\begin{equation}\label{encl2}
(v_+f)(Z)= f(\pi(Z)).
\end{equation}

\noindent
{\bf Remark}.
As it was mentioned, the function $b$ is not single valued
but $|b(z)|$ is a single valued function. We define the character $\mu\in
\Gamma^*$ by
$$
b(\gamma z)=\mu(\gamma) b(z).
$$
Let $\gamma_j$ be  the contour, that starts at infinity (or any
other real point bigger than $1$), go in the upper half--plane to
the gap $(a_j,b_{j})$ and then go back in the lower half--plane to
the initial point. Assuming that
$b_0<...<a_j<b_j<a_{j+1}<...<a_0$, we have
$\mu(\gamma_j)=e^{-{2\pi i}\frac{g+1-j}{g+1}}$, equivalently
$\omega([b_j,a_{0}])=\frac {g+1-j} {g+1}$.
Note that the system of the above contours $\gamma_j$
is a generator of the free group $\Gamma^*(E)$. In other words
a character $\alpha$ is uniquely defined by the vector
$$
\begin{bmatrix}
\alpha(\gamma_1),&\alpha(\gamma_2),&\dots,&\alpha(\gamma_g)
\end{bmatrix}\in \TT^g.
$$
This sets an one--to--one correspondence between
$\Gamma^*(E)$ and $\TT^g$.

Recall the key role of the reproducing kernels
$k^{\alpha}$ in our  construction. In this particular case
they especially
well studied
\cite{Fay}. First of all, they have analytic continuation (as
multivalued functions) on the whole
$X$, so we can write
$k^{\alpha}(Z)$.

\begin{theorem} For every $\alpha\in \Gamma^*$
  the reproducing kernel $k^{\alpha}(Z)$ has on $X$ exactly
$g$ simple poles that do not depend on
$\alpha$ and
$g$ simple zeros. The divisor $D_+= \sum_{j=1}^g Z_j$
of zeros
\begin{equation}\label{4.9}
k^{\alpha}(Z_j)=0
\end{equation}
with the divisor of poles form the divisor
\begin{equation}\label{4.14s}
\div (k^{\alpha})= D_+- D_C
\end{equation}
that belongs to $D(E)$, moreover \eqref{4.14s} sets an  one--to--one
correspondence
between $D(E)$ and $\Gamma^*(E)$.
\end{theorem}
The functions $k^{\alpha}$ possess different representations, in
particular, in terms of theta--functions \cite{Mam}, and the map
$D\mapsto \alpha$ can be written explicitly in terms
of abelian integrals (the Abel map).

\bigskip

\noindent
{\bf Summary}. The three objects $J(E)$, $D(E)$ and $\Gamma^*(E)$
are equivalent. Both maps  $\Gamma^*(E)\to D(E)$
and $\Gamma^*(E)\to J(E)$
can be defined in terms of the reproducing kernels of the spaces
$H^2(X_+,\alpha)$,
$\alpha\in \Gamma^*(E)$. The first one is given by
\eqref{4.14s}. It associates to the given $k^{\alpha}(Z)$ the sets
of its zeros and poles (the poles are fixed and the zeros vary with $\alpha$).
The
matrix $J(\alpha)\in J(E)$ is defined as the matrix
of the multiplication
operator by $z(Z)$
\begin{equation}\label{jm}
z(Z)e^{\alpha}_s(Z)=p^\alpha_s e^{\alpha}_{s-1}(Z)+
q^\alpha_s e^{\alpha}_{s}(Z)+
p^\alpha_{s+1} e^{\alpha}_{s+1}(Z),\ Z\in X, \ s\in\ZZ,
\end{equation}
with respect to the basis
\begin{equation}\label{basisjm}
e^{\alpha}_{s}(Z)=b^s(Z) K^{\alpha\mu^{-s}}(Z).
\end{equation}
\bigskip

It's really easy to see that $J(\alpha)$
is periodic: just recall that $b^{g+1}$ is single valued, that is,
$\mu^{g+1}=1$, and therefore the spaces $H^2(X_+,\alpha)$
and $H^2(X_+,\alpha\mu^{-(g+1)})$ (and their reproducing kernels)
coincide. \bigskip

Now we can go back to the Renormalization Equation. Note that
$\pi$ acts naturally on $\Gamma(E)$:
$$
\pi\gamma =\{\pi(Z),\ Z\in\gamma\}\in \Gamma(\tilde E),
\quad \text{for}\ \gamma\in \Gamma(E).
$$
The map $\pi^*: \Gamma^*(\tilde E)\to \Gamma^*(E)$ is defined by
duality:
\begin{equation}\label{upch}
(\pi^* \tilde\alpha)(\gamma)= \tilde \alpha (\pi\gamma).
\end{equation}

\begin{theorem} \label{th4.7}
Let $T$, $T^{-1}:[-1,1]\to [-1,1]$,
be an expanding polynomial.
Let $\tilde J$ be a periodic Jacobi matrix with spectrum $\tilde E \subset
[-1,1]$, and therefore there exists a polynomial $\tilde U$ such that $\tilde
E=\tilde U^{-1} [-1,1]$ and a character $\tilde\alpha\in \Gamma^*(\tilde E)$
such that $\tilde J= J(\tilde \alpha)$. Then
\begin{equation}\label{solu1}
J:=J(\pi^*\alpha)=\mathcal J_T(\tilde J)
\end{equation}
is the periodic Jacobi matrix with spectrum $E= U^{-1} [-1,1]$, $U:=\tilde
U\circ T$, that satisfies the Renormalization Equation
\begin{equation}\label{t01}
V^*(z-J)^{-1}V=\frac{T'(z)} d (T(z)-\tilde J)^{-1},
\end{equation}
where the isometry matrix $V$ is defined by $V|k\rangle:=|kd\rangle$.
\end{theorem}

\begin{remark}
Let us mention that the Renormalization Equation can be rewritten
equivalently in the  form of polynomials equations, as it should be
since  we have
the map from one algebraic variety, $\Jac(\tilde X)$, in the another
one,  $\Jac(X)$. Equation \eqref{t01} is equivalent to, see
\cite{PVYU},
\begin{equation}
\label{re.1}
V^*T(J)=\tilde JV^*,
\end{equation}
\begin{equation}
\label{re.2}
V^*\frac{T(z)-T(J)}{z-J}V=T'(z)/d.
\end{equation}
\end{remark}

\begin{proof}[Proof of Theorem \ref{th4.7}]
First we note, that for the operator multiplication
by $z(Z)$ in $L^2(\partial X_+)$, the operator multiplication
by $\tilde z(\tilde Z) $ in $L^2(\partial \tilde X_+)$,
the spectral parameter $z_0$ and the isometry
\begin{equation*}
(v f)(Z)=f(\pi(Z)),\quad v:L^2(\partial \tilde X_+)\to
L^2(\partial X_+),
\end{equation*}
we have
\begin{equation}\label{4.20s2}
\int_{\partial X_+} \frac 1{z_0-z(Z)} |(v f)(Z)|^2\, d\omega=
\int_{\partial \tilde X_+}\left(\frac  1 d
\sum_{\pi(Z)=\tilde Z}
\frac 1{z_0-z(Z)}\right) | f(\tilde Z)|^2\, d\tilde\omega.
\end{equation}
It is evident, that
\begin{equation*}
\frac  1 d
\sum_{T(y)=x}
\frac 1{z_0-y}=\frac{T'(z_0)/d}{T(z_0)-x}.
\end{equation*}
Thus
\begin{equation}\label{4.21s2}
v^*(z_0-z(Z))^{-1}v=(T'(z_0)/d)(T(z_0)-\tilde z(\tilde Z))^{-1}.
\end{equation}

It remains to show that $\pi$ transforms the basis vector
$
\tilde e_n^{\tilde\alpha}=
\tilde b^n K^{\tilde\alpha\tilde\mu^{-n}}
$
into
$$
  e_{nd}^{\pi^*\alpha}=
  b^{nd} K^{(\pi^*\tilde\alpha)\mu^{-n d}}=
(\tilde b^n\circ\pi) K^{\pi^*(\tilde\alpha\tilde\mu^{-n})}.
$$
Or, what is the same, that
$
  K^{\tilde\alpha}\circ\pi=K^{\pi^*\alpha}
$ for all $\tilde \alpha\in\Gamma^*(\tilde E)$.
Note that both functions are of norm one in the
same space $H^2(X_+,\pi^*\tilde \alpha)$, in particular, they have the
same character of automorphy $\pi^*\tilde \alpha\in \Gamma^*(E)$. Note,
finally, that the divisor
$$
\div (k^{\tilde\alpha}\circ\pi)= \pi^{-1} (\tilde D_+)-
\pi^{-1} (\tilde D_C),
$$
where
$\div (k^{\tilde\alpha})= \tilde D_+-\tilde D_C$,
belongs to $D(E)$, therefore $k^{\tilde\alpha}\circ\pi$ is the reproducing
kernel and the theorem is proved.
\end{proof}

Probably it would be better to call \eqref{t01} the Renormalization
Identity in the above theorem. The idea is that one can try
to define $J$ as the {\it solution} of \eqref{t01}
with the given $\tilde J$. Indeed, in the case of one--sided
matrices such equation has the unique solution.
Now we demonstrate that in two sided case
for the given periodic
$\tilde J$ we can find $2^{d-1}$ solutions.

To find all this solutions of
\eqref{t01} let us look  a bit more carefully at the above proof.
Note that the same identity \eqref{4.20s2} holds for any
isometry $v$ of the form
$$
vf=v_\theta f=\theta (f\circ\pi),
$$
where $\theta$ is a unimodular ($|\theta|=1$) function
on $\partial X_+$.

Concerning the second part of the proof, let us mention that
the set of critical points of $U$ splits in two sets:
$$
\{c_U\}=T^{-1}\{c_{\tilde U}\}\cup \{c_T\}.
$$
Correspondingly,
$$
\sum (C_U)_j=\sum_{k}
\sum_{\pi (C_{ U})_{k,j}=(C_{\tilde U})_k}(C_U)_{k,j}+\sum (C_T)_j,
$$
and the divisor of $k^{\tilde\alpha}\circ\pi$
consists of two parts, that one that depends on $\tilde \alpha$
$$
\pi^{-1}(\tilde D),
$$
and that part that corresponds to the critical points
of the polynomial $T$
$$
\{(C_T)_j\}_{j=1}^{d-1},
$$
since
$$
D=\div(k^{\tilde\alpha}\circ\pi)=
\pi^{-1}(\tilde D)+\sum_{j=1}^{d-1}(C_T)_j-
\pi^{-1}(\tilde D_C)-\sum_{j=1}^{d-1}(C_T)_j.
$$
Thus we can fix an arbitrary system of points $\{Z_{c,j}\}_{j=1}^{d-1}$ such
that
$z(Z_{c,j})$ belongs to the same gap in the spectrum $E$ as the critical point
$(c_T)_j$. If $\theta$ is the canonical product on $X$ with the divisor
$$
\div(\theta)=\sum_{j=1}^{d-1} Z_{c,j}-\sum_{j=1}^{d-1} (C_T)_{j},
$$
then $\theta k^{\tilde\alpha}\circ\pi$ is the reproducing kernel
simultaneously for all $\tilde \alpha\in\Gamma^*(\tilde E)$.
But to make $\theta$ unimodular
(zeros and poles are symmetric)
our choice is restricted just to
$Z_{c,j}=(C_T)_{j}$ or $Z_{c,j}=\tau (C_T)_{j}$. Note that
$\tau (C_T)_{j}- (C_T)_{j}$ is the devisor of the Complex Green function
$b_{(c_T)_j}$. In this way we arrive at

\begin{theorem}\label{th5.8} For an expanding polynomial $T$, and a
periodic Jacobi matrix
$\tilde J= J(\tilde \alpha)$, $\tilde\alpha\in\Gamma^*(\tilde E)$ as in Theorem
\ref{th4.7} there exist $2^{d-1}$ solutions of the Renormalization Equation
\eqref{t01}. Denote by $\mu_{(c_T)_j}$ the character generated by the Green's
function $b_{(c_T)_j}$, $ b_{(c_T)_j}\circ\gamma= \mu_{(c_T)_j}(\gamma)
b_{(c_T)_j}. $ Then these solutions are of the form
\begin{equation}\label{solu2}
J:=J\left(
\eta_\delta
\pi^*\tilde\alpha\right),
\quad  \eta_\delta:=\prod^{d-1}_{j=1}\mu_{(c_T)_j}^{\frac 1 2
(1+\delta_{(c_T)_j})},
\end{equation}
as before
\begin{equation*}
     \delta=\{\delta_{(c_T)_j}\},\quad \delta_{(c_T)_j}=\pm 1.
\end{equation*}
\end{theorem}

\begin{proof}
We define the isometry
$$
(vf)(Z)=\left(\prod^{d-1}_{j=1}
b_{(c_T)_j}^{\frac 1 2(1+\delta_{(c_T)_j})}(Z)\right)
f(\pi(Z))
$$
and then repeat the arguments of the proof of Theorem \ref{th4.7}.
\end{proof}

\section {The Renormalization Equation for two--sided
Jacobi matrices (general case)}

In this section
we assume that
   $$
   T(z)=z^d-qdz^{d-1}+...
   $$
   is a {\it monic}
  expanding polynomial. Under this normalization
  $T^{-1}:[-\xi,\xi]\to [-\xi,\xi]$ for a certain $\xi>0$.
It was proved in \cite{PVYU} that the Renormalization Equation
has $2^{d-1}$ solutions for every two
sided Jacobi matrix $\tilde J$ with the spectrum
in $[-\xi,\xi]$, not only for periodic one.
Moreover,
they are the only possible solutions.
For the reader convenience we formulate these theorems here.

By $l^2_{\pm}(s)$ we denote the spaces
  which are formed by
$\{|s+k\rangle\}$ with $k\le 0$ and $k\ge 0$ respectively, that
is $l^2(\ZZ)=l^2_-(s)\oplus l^2_+(s+1)$. Correspondingly to these
decompositions we set
$\tilde J_{\pm}(s)=P_{l^2_{\pm}(s)}\tilde
J|l^2_{\pm}(s)$.
Recall that  a (finite or infinite)
one--sided Jacobi matrix is uniquely determined by its
so called resolvent function
\begin{equation}\label{4aug22}
\tilde r_{\pm}(z,s)=\langle s|(\tilde J_{\pm}(s)-z)^{-1}
|s\rangle.
\end{equation}
%


This set of solutions we parametrize
by a collections of vectors
\begin{equation}\label{delta}
\delta:=\{\delta_c\}_c,
\end{equation}
where each component $\delta_c$ can be chosen as plus or
minus one.

\begin{theorem}\label{exs} Fix a vector $\delta$ of the form
\eqref{delta}.
For a given $\tilde J$ with the spectrum on $[-\xi,\xi]$ define the
Jacobi matrix $J$
according to the following algorithm:

For $s\in \ZZ$ we put
\begin{equation}\label{4m}
\frac 1{T^{(s)}(c)}=-\tilde r_-(T(c),s),
\quad\text{if}\  \delta_c=-1,
\end{equation}
  and
\begin{equation}\label{4p}
{T^{(s)}(c)}=-\tilde p_{s+1}^2\tilde r_+(T(c),s+1),
\quad\text{if}\  \delta_c=1,
\end{equation}
  where the functions $\tilde r_\pm(z,s)$
are defined by \eqref{4aug22}. Then define the monic polynomial
$T^{(s)}(z)$ of degree $d$ by the interpolation formula
\begin{equation}\label{4.2}
T^{(s)}(z)=(z-q)T'(z)/d+
\sum_{c:T'(c)=0}\frac{T'(z)}{(z-c)T''(c)}T^{(s)}(c).
\end{equation}
Define the block

\begin{equation}\label{4aug23}
J^{(s)}=
\begin{bmatrix}
q_{sd}&p_{sd+1}& & & \\
p_{sd+1}&q_{sd+1}&p_{sd+2}& & \\
     &\ddots&\ddots&\ddots & \\
& & p_{sd+d-2}&q_{sd+d-2}&p_{sd+d-1} \\
&  & &p_{sd+d-1}&q_{sd+d-1}
\end{bmatrix}
\end{equation}
by its resolvent function
\begin{equation}\label{3*}
\left<0\left|(z-J^{(s)})^{-1} \right|0\right>=
\frac{T'(z)/d}{T^{(s)}(z)},
\end{equation}
where $T^{(s)}(z)$ is a monic polynomial of degree $d$.
 Finally define the entry $p_{sd+d}$ by
$
p_{sd+1}...p_{sd+d}=\tilde p_{s+1}$

We claim that the matrix $J=J(\delta, \tilde J)$, combined with such
blocks and entries over all $s$, satisfies \eqref{t01}.
\end{theorem}

\begin{theorem}\label{uni}
Theorem \ref{exs} describes the whole set of solutions of the
Renormalization Equation.
\end{theorem}

In \cite{PVYU} we concentrate only on one of the
solutions of
\eqref{t01}, namely that one that  related to the vector
$$
\delta_-=\{-1,\dots ,-1\},
$$
that is all $T^{(s)}(c)$ are defined by \eqref{4m}.
Note, $J(\delta_-,\tilde J)=\mathcal J_T(\tilde J)$
for a periodic $\tilde J$.
Precisely for this
solution
  we proved
   (main Theorem 1.1 in \cite{PVYU}):
  \begin{theorem} \label{mainth}
Let $\tilde J$ be a Jacobi matrix
with the spectrum on $[-\xi,\xi]$. Then the  Renormalization
Equation \eqref{t01} has a
solution $J=J(\delta_-,\tilde J)$
with the spectrum on
$T^{-1}([-\xi,\xi])$.
Moreover,
  if
$\min_i |t_i|\ge 10\xi$ then
\begin{equation}\label{contr1}
  \Vert J(\delta_-,\tilde J_1)-J(\delta_-,\tilde J_2)\Vert
  \le \kappa\Vert \tilde J_1-\tilde J_2\Vert
  \end{equation}
with an absolute constant $\kappa<1$.
\end{theorem}
Let us emphasize the especial role of the position of
the critical values: the transform
$\mathcal J_T(\tilde J)$ is a contraction as soon as critical values are
distant sufficiently
far from the spectrum.

  In this work we add certain remarks concerning other
  solutions of \eqref{t01}.


\subsection{The duality $\delta\mapsto-\delta$}
At least for one more solution of the Renormalization Equation
is a contraction.

\begin{theorem}\label{th6.1}
The dual solution of the Renormalization Equation $J(\tilde J,-\delta)$,
  possesses the contractibility
property simultaneously with  $J(\tilde J,\delta)$.
\end{theorem}
It deals with the
following universal involution acting on  Jacobi matrices
\begin{equation}  \label{is6.1}
J\to J_\tau:=U_\tau J U_\tau, \quad\text{where}\ U_\tau
|l\rangle=|1-l\rangle.
\end{equation}
Obviously $VU_\tau=U_\tau S^{1-d}V$. Thus, having $J$ as a solution of
the renormalization equation corresponding to $\tilde J$ we have
simultaneously that $S^{d-1}J_\tau S^{1-d}$  solves the equation with
the initial $\tilde J_\tau$.  The following lemma describes which branch
corresponds
to which in this case.

\begin{lemma}
Let $J=J(\tilde J,\delta)$ then
\begin{equation}\label{is6.2}
   S^{d-1}J_\tau S^{1-d}=J(\tilde J_\tau,-\delta).
\end{equation}
\end{lemma}
\begin{proof}
We give a proof using the language of Sect. 5, so formally we prove the claim
only for periodic matrices.

Note that the involution
\eqref{is6.1} is strongly related to the standard involution $\tau$
\eqref{invol} on $X$.
Indeed, the function $K(\tau Z, \alpha)$ has the divisor
$$
\tau D_+ -\tau D_C=(\tau D_+ - D_C)-(\tau D_C-D_C),
$$
that is,
\begin{equation*}
K(\tau Z, \alpha)=\frac{K( Z, \beta)}{b_{c_1}(Z)\dots b_{c_{g}}(Z)},
\end{equation*}
and $\beta=\nu\alpha^{-1}$, where $\nu=\mu_{c_1} \dots   \mu_{c_g} $.
Due to this remark and the property $z(\tau Z)=z(Z)$ we have
\begin{equation}\label{is6.3}
(J(\alpha))_\tau=J(\nu\mu\alpha^{-1}).
\end{equation}

Now we apply \eqref{is6.3} to prove \eqref{is6.2}.  Let $\tilde
J_\tau=J(\tilde\alpha)$ with $\tilde\alpha\in \Gamma^*(\tilde X_+)$. Or,
in other words, $\tilde J=J(\tilde\mu\tilde\nu\tilde\alpha^{-1})$.
Then by \eqref{solu2}
\begin{equation*}
     J(\tilde
     J,\delta)=J(\eta_\delta\pi^*(\tilde\mu\tilde\nu\tilde\alpha^{-1})),
     \quad  \eta_\delta:=\prod^{d-1}_{j=1}
\mu_{(c_T)_j}^{\frac 1 2 (1+\delta_{(c_T)_j})}.
\end{equation*}
But $\pi^*\tilde\mu=\mu^d$ and $\pi^*(\tilde\nu)=\nu\eta_{\delta_+}^{-1}$
(just to look at the characters of the corresponding Blaschke products).
Thus, having in mind that $\eta_{\delta}  \eta_{-\delta}= \eta_{\delta_+}$,
we obtain
\begin{equation*}
    J(\tilde
    J,\delta)=J(\mu^d\nu\eta^{-1}_{-\delta}\pi^*(\tilde\alpha^{-1})).
\end{equation*}
Using again \eqref{is6.3} we get
\begin{equation*}
    (J(\tilde
    J,\delta))_\tau=J(\mu^{1-d}\eta_{-\delta}\pi^*(\tilde\alpha))
    =S^{1-d}J(\eta_{-\delta}\pi^*(\tilde\alpha))S^{d-1},
\end{equation*}
and the lemma and Theorem \ref{th6.1} are proved.
\end{proof}

Having two different contractive branches of solutions of the
renormalization equation, following \cite{Kn}, to an arbitrary sequence

\begin{equation*}
\epsilon=\{\epsilon_0,\epsilon_1,\epsilon_2\dots\},
\quad \epsilon_j=\delta_{\pm}.
\end{equation*}
we can associate a limit periodic
matrix $J$ with the spectrum on $\Julia(T)$. For a fixed sufficiently
hyperbolic
polynomial $T$, we define $J$ as the limit of
\begin{equation}\label{is6.4}
J_n:=J(
\eta_{\epsilon_0} \pi^*  \eta_{\epsilon_1} \dots   \pi^*\eta_{\epsilon_{n-1}}).
\end{equation}

\subsection{Other solutions of the Renormalization Equation and the
Ruelle operators}
Iterating the renormalization transform
$$
 J_{n+1}:=J(J_n,\delta_-),
$$
 we obtain almost periodic operator $J=\lim J_n$
with the spectrum on the Julia set of the sufficiently expanding
polynomial
$T$ (the key point here, of course  property \eqref{contr1}).  This
operator  possesses the following structure \cite{PVYU}: it is the direct
sum of two (one--sided) Jacobi matrices $J_{\pm}$. Moreover, the
spectral measure of the operator $J_+$ is  the balanced measure $\mu$
on the Julia set. That is it is the eigen--measure of the Ruelle
operator $\mathcal L^*$, where
$$
(\mathcal L f)(x)=\frac 1{d}\sum_{T
n(y)=x}f(y).
$$
The spectral measure $\nu$ of $J_-$ is the so called
 Bowen--Ruelle measure. It is
 the eigen--measure
for $\mathcal L_2^*$,
$$
(\mathcal L_2 f)(x)=\sum_{T(y)=x}\frac{f(y)}{T'(x)^2},
$$

We conjecture that actually all branches of
solutions of the renormalization equation are contractions for
sufficiently hyperbolic $T$. At least the previous subsection looks as
a  quite strong indication in this direction: considering,
instead of initial $T$,  $T^2=T\circ T$ or its bigger powers, we
get, as in \eqref{is6.4}, several $\delta$'s,
  $ \eta_{\delta} =\eta_{\epsilon_0} \pi^*  \eta_{\epsilon_1} \dots
  \pi^*\eta_{\epsilon_{n-1}}$,
  possessing the contractibility property with respect to the polynomial $T^n$
and different
from $\delta_{\pm}$ (related to $(\pi^*)^n$).

  Similarly to the above statement
we formulate

\smallskip

\noindent
{\it Conjecture}.
     Let $T(z)$ be an expanding polynomial. For a given $\delta$
     we factorize
     $T'(z)=A_1(z) A_2(z)$ putting in the first factor
     all critical points related to $\delta_c=1$. Denote by
$\sigma_{1,2}$,  the
     (nonnegative)
     eigen--measures, corresponding to the Ruelle operators
     \begin{equation}\label{}
     (\mathcal L_{A_i}f)(x)=\sum_{T(y)=x}\frac {f(y)}{A_i(y)^2},
\end{equation}
i.e.,  $\mathcal L^*_{A_i}\sigma_i=\rho_i  \sigma_i $. Finally let
$J_{1,2}$ be the one-sided Jacobi matrices associated with
$\sigma_{1,2}$.
We conjecture that the iterations
$
 J_{n+1}:=J(J_n,\delta),
$
converges to  the
block matrix $J=J_-\oplus J_+$ with $J_-=J_1$ and $J_+=J_2$.
In particular this means that all such operators are limit
periodic.
\bigskip

Note that  the conjecture holds true also,
   say, for $T^2(z)$ and $A_1(z)=T'(z)$, $A_2= T'(T(z))$.

\subsection{Shift transformations with the Lipschitz property}
We say that the direction $\eta\in \Gamma^*$ has   the Lipschitz property
with a constant $C(\eta)$ if for all   $\alpha, \beta\in \Gamma^*$
\begin{equation}\label{}
     \Vert J(\eta\alpha)-J(\eta\beta)\Vert\le   C(\eta)
     \Vert J(\alpha)-J(\beta)\Vert.
\end{equation}
Then, one can get the contractibility of the map $\eta\pi^*$ in two steps:
\begin{equation}
\begin{split}
     \Vert J(\eta\pi^*\tilde\alpha)-J(\eta\pi^*\tilde\beta)\Vert &\le   C(\eta)
     \Vert J(\pi^*\tilde\alpha)-J(\pi^*\tilde\beta)\Vert   \\
     &\le C(\eta)\kappa
     \Vert J(\tilde\alpha)-J(\tilde\beta)\Vert.
     \end{split}
  \end{equation}

Note, that in fact the situation is a bit more involved because we
should be able to
compare Jacobi matrices with different spectral sets, for example, when
$E_i=T^{-1}\tilde E_i$, $\tilde E_1\not=\tilde E_2$. But
we just want to
indicate the general idea. Note, in particular, that for directions
$\eta_{\delta}$ of the form \eqref{solu2} such a comparison is
possible. Of course, for our goal the constant $C(\eta)$ should be
uniformly bounded when we increase the level of sufficient
hyperbolicity of $T$ making $\kappa$ smaller.

However the key point of this remark (this way of proof) is that,
actually, {\it we do not need to constrain ourselves   by the
form of the vector $\eta$}. Combining a  ``Lipschitz" shift by $\eta$ (the
direction is restricted just by this property) with a sufficiently
contractive pull--back $\pi^*$ we arrive at an iterative process that produces
a limit periodic Jacobi matrix with the spectrum on the same
$\Julia(T)$. In the next
subsection we give examples of directions with the required property, see
Corollary \ref{cordbr}.

We do not
have a proof of the Lipschitz property of
$\eta_\delta$'s, but there is a good chance to generalize the result
of the next
subsection in a way that at least some of the directions
$\eta_\delta$ will be also
available.

Finally, we would be very interested to know, whether there is in
general a relation
between the form of the ``weight" vector $\eta$ and the corresponding
weights of the
Ruelle operators (if any exists).

\subsection{Quadratic polynomials and
the Lipschitz property of the Darboux transform}
Consider the simplest special case $T(z)=\rho(z^2-1)+1$, $\rho>2$.
Note that  the spectral set
$E=T^{-1}\tilde E$
is symmetric, moreover
the matrix related to $H^2(\pi^*\tilde\alpha)$ has zero main diagonal
(as well as a one--sided matrix related to a symmetric measure).
Now we introduce a decomposition of $H^2(\pi^*\tilde\alpha)$
which is
very similar to the standard decomposition into even and odd functions.

We define the two--dimensional vector--function representation of $f\in
H^2(\pi^*\tilde\alpha)$
\begin{equation}\label{30s36}
f\mapsto\frac 1{\sqrt{2}}
\begin{bmatrix} f(Z_1(\tilde Z)) \\ f(Z_2(\tilde Z))
\end{bmatrix}\mapsto
\begin{bmatrix} g_1(\tilde Z)\\ g_2(\tilde Z)
\end{bmatrix},
\end{equation}
where
\begin{equation*}
\begin{bmatrix} g_1(\tilde Z)\\ g_2(\tilde Z)
\end{bmatrix}=
\frac 1 2
\begin{bmatrix} f(Z_1(\tilde Z)) +f(Z_2(\tilde Z))\\
f(Z_1(\tilde Z))-f(Z_2(\tilde Z))
\end{bmatrix},
\end{equation*}
the first component, in a sense, is even and the second is odd.
To be more precise, let us describe analytical properties of this
object in $\partial \tilde X_+$.

Note that
due to
$$
\int_{\partial X_+}|f|^2 d\omega=
\int_{\partial \tilde X_+}\frac 1 2\sum_{\pi(Z)=\tilde Z}|f|^2 d
\tilde\omega
$$
metrically it is of $L^2$ with respect to $\tilde\omega$,
moreover the transformation is norm--preserved.

It is evident that the function
$g_1$
belongs to $H^2(\tilde X_+, \tilde\alpha)$. Consider the second function.
Note that the critical points of $T$ are zero and infinity.
For a small circle $\gamma$ around the point
$T(0)=-\rho+1$ we have
$g_2\circ\gamma=-g_2$
and the same property for a contour $\gamma$ that surrounds infinity. Let us
introduce
$$
\Delta^2:= \tilde b_{T(0)} \tilde b.
$$
Note that for the above contours we have $\Delta\circ\gamma=- \Delta$.
We are going to represent $g_2$ in the form $g_2=\Delta\hat g_2$ and to claim
that
$\hat g_2$ has nice automorphic properties in $\tilde X_+$. Let us note that
\begin{equation*}
\tilde b\frac{\tilde z-T(0)}{\tilde b_{T(0)}}=
\tilde b^2\frac{\tilde z-T(0)}{\Delta^2}
\end{equation*}
is an outer function in the domain $\bar{\CC}\setminus \tilde E\simeq
\tilde X_+$.
So, the square root of this function is well defined. We put
\begin{equation}\label{sqroot}
\tilde b \phi:=\sqrt{
\tilde b^2\frac{\tilde z-T(0)}{\rho \Delta^2}}
\end{equation}
and denote by $\tilde \eta$ the character generated by $ \phi$,
$ \phi\circ\gamma= \eta(\gamma) \phi$. Thus
\eqref{sqroot} reduces the ramification of the function
$\Delta$ to the function $\phi$, which is well defined in the domain,
and to the
elementary function $\sqrt{\tilde z-T(0)}$.

\begin{theorem}
The transformation $f\mapsto g_1\oplus \hat g_2$ given by
\eqref{30s36} is a unitary map from $H^2(\pi^*\tilde\alpha)$
to $H^2(\tilde\alpha)\oplus H^2(\tilde\alpha\tilde\eta)$.
Moreover with respect to this representation
\begin{equation}\label{12o40}
z f\mapsto
\begin{bmatrix}0& \bar\phi\\
\phi &0
\end{bmatrix}\begin{bmatrix}
g_1\\ \hat g_2
\end{bmatrix}
\end{equation}
and
$$
v_+f\mapsto f\oplus 0,\quad f\in H^2(\tilde\alpha),
$$
where the isometry $v_+:H^2(\tilde\alpha)\to H^2(\pi^*\tilde\alpha)$
is defined by \eqref{encl2}.
\end{theorem}

\begin{proof}
By the definition of $\Delta$ we have
\begin{equation}\label{12o38}
g_2=\Delta \hat g_2,\
\text{where}\ \hat g_2\in H^2(\tilde\alpha\tilde\eta).
\end{equation}
Further,
since
$$
  z_{1,2}=\pm\sqrt{\frac{\tilde z-T(0)}{\rho}},
$$
we have, say for the second component,
\begin{equation}\label{12sept38}
\frac{1}{\Delta}\frac{(zf)(Z(\tilde Z_1))- (zf)(Z(\tilde Z_2))}{2}
=\sqrt{\frac{\tilde z-T(0)}{\rho\Delta^2}}
\frac{f(Z(\tilde Z_1))+ f(Z(\tilde Z_2))}{2}=
\phi g_1.
\end{equation}
Since on the boundary of the domain
$$
\phi^2\Delta^2=\frac {\tilde z-T(0)}{\rho}=|\phi|^2
$$
(the second expression is positive on $\partial \tilde X_+$)
we have
\begin{equation}\label{12o41}
\phi\Delta^2=\overline{\phi} \quad {\rm on}\ \tilde E.
\end{equation}
Using this relation, similarly to \eqref{12sept38}, we prove the identity
of the first components in \eqref{12o40}.

\end{proof}

\begin{theorem}
The multiplication operator
$\phi: L^2(\partial \tilde X_+)\to L^2(\partial \tilde X_+)$
with respect to the basis systems \eqref{basisjm} related to
$\tilde\alpha$ and
$\tilde \eta\tilde \alpha$, respectively, is a two diagonal matrix $\Phi$.
Moreover,
\begin{equation}\label{12o42}
\Phi^*\Phi= \frac{J(\tilde \alpha)-T(0)}{\rho},
\quad
\Phi\Phi^*= \frac{J(\tilde\eta\tilde \alpha)-T(0)}{\rho}.
\end{equation}
In other words, the transformation
$J(\tilde \alpha)\mapsto J(\tilde\eta\tilde \alpha)$ is the Darboux
transform.
\end{theorem}

\begin{proof} First of all $\phi$ is a character--automorphic function with the
character $\tilde\eta$ with a unique pole at infinity ($\tilde b\phi$ is an
outer function). Therefore the multiplication operator acts from $\tilde
bH^2(\tilde\alpha \tilde\mu^{-1})$ to $H^2(\tilde\eta\tilde\alpha)$. Therefore,
the operator $\Phi$ has only one non--trivial diagonal above the main diagonal.
The adjoint operator has the symbol $\overline \phi$. According to
\eqref{12o41}
it has holomorphic continuation from the boundary inside the domain. Thus
$\Phi^*$ is a lower triangular matrix. Combining these two facts we get that
$\Phi$ has only two non--trivial diagonals. Then, just comparing symbols of
operators on the left and right parts of \eqref{12o42}, we prove these
identities.
\end{proof}

\begin{corollary}\label{cordbr}
    Let $\tilde J_{1,2}$ be periodic Jacobi matrices with the
spectrum on $[-1,1]$. Let
    $\Drb(\tilde J_{1,2},\rho)$ be their Darboux transforms. Then
\begin{equation}\label{drx}
\Vert\Drb(\tilde J_{1},\rho)-\Drb(\tilde J_{2},\rho)\Vert
\le C(\rho)
\Vert\tilde J_{1}-\tilde J_{2}\Vert.
\end{equation}
\end{corollary}
\begin{proof}
For the given $\tilde J_{1,2}$ we define $J_{1,2}$ via the quadratic polynomial
$T(z)=\rho(z^2-1)+1$. Being decomposed into even and odd indexed subspaces they
are of the form
\begin{equation}\label{drx2}
J_{1,2}=\begin{bmatrix}
0&\Phi^*_{1,2}\\
\Phi_{1,2}& 0
\end{bmatrix}.
\end{equation}
Due to the main theorem, that gives the uniform estimate for
$\Vert J_1-J_2\Vert$, we have
\begin{equation}\label{drx}
\Vert\Phi_{1}-\Phi_{2}\Vert
\le \kappa(\rho)
\Vert\tilde J_{1}-\tilde J_{2}\Vert,
\end{equation}
with $\kappa(\rho)=\frac{C}{\rho-2}$, $C$ is an absolute constant.
Using \eqref{12o42} we get \eqref{drx} with
$C(\rho)=\frac{2\rho C}{\rho-2}$.
\end{proof}

\section{double covering $\pi(v)=\tau v-1/v$}

\subsection{Over the
simply--connected domain}
The simplest expanding double coverings are $T_1=v^2-\lambda$,
$\lambda>2$, and $T_2=\tau v-1/v$, $\tau>1$. Denote by $\xi_{1,2}$
the fixed points $T_{1,2}(\xi_{1,2})=\xi_{1,2}$, $\xi_{1,2}>0$.
To start with we give a complete description of finite difference
operators related to these coverings over the
simply--connected domain $\overline{\CC}\setminus E$, where
$E=[-\xi_{1,2},\xi_{1,2}]$.

More precisely, we start with
a  Jacobi matrix $J_0$ with constant coefficients.
Under the normalization $\sigma(J_0)=E$ we have
$J=\frac{\xi_{1,2}} 2(S+S^{-1})$. That is the symbols of this
operator $(z_*,b_*)$ are related by
\begin{equation}\label{20s17}
z_*=\frac{\xi_{1,2}} 2\left(\frac 1{b_*}+b_*\right),
\end{equation}
($b_*$ is the Green function of the domain $\overline{\CC}\setminus E$).
For $\pi=T_1$ or $\pi=T_2$ we consider
the open Riemann surfaces $Y_c\setminus F$
with $Y_c\simeq\overline{\CC}$, $F=\pi^{-1}(E)$
and describe operators with the  symbols $(z^*,b^*)$:
\begin{equation}\label{20s16}
z^*=z_*\circ \pi\quad\text{and}\quad
(b^*)^2=b_*\circ \pi.
\end{equation}

The main difference between these two cases is that in the first
one we have only one infinity on $Y_c$ ($\infty\in \overline{\CC}$)
and in the second case there are two infinities:
$0,\infty\in\overline{\CC}$.
Correspondingly an intrinsic basis contains the reproducing kernels
related only to one fixed point in the first case and to
two specific points in the second case. As result the multiplication
operator by $v$ with respect to this basis is a Jacobi matrix
in the first case and a five diagonal matrix (of a special form,
see Lemma 2.1) in the second.

\subsubsection{$\pi=v^2-\lambda$}
Due to \eqref {20s17}, \eqref{20s16} we have
\begin{equation}\label{20s18}
v^2-\lambda=\frac{\xi_1} 2\left(\frac 1{(b^*)^2}+(b^*)^2\right).
\end{equation}
Since $v$ is the symbol of a Jacobi matrix,
\begin{equation}\label{20s19}
v\sim S\Lambda_1 +\Lambda_0 +\Lambda_1 S^{-1},
\end{equation}
where $\Lambda_{0,1}$ are diagonal of period two matrices,
and $(b^*)^2\sim S^2$, we have from \eqref{20s18}
\begin{equation}\label{20s20}
\begin{split}
&\lambda_0^{(1)}\lambda_1^{(1)}=\xi_1/2,\\
&\lambda_0^{(0)}+\lambda_1^{(0)}= 0,\\
&(\lambda_0^{(0)})^2+
(\lambda_0^{(1)})^2+(\lambda_1^{(1)})^2
=\lambda.
\end{split}
\end{equation}

\subsubsection{$\pi=\tau v-1/v$} Repeating arguments \eqref{20s16},
\eqref{20s17} we get instead of \eqref{20s18}
\begin{equation*}
\tau v-1/v=\frac{\xi_2} 2\left(\frac 1{(b^*)^2}+(b^*)^2\right).
\end{equation*}
or
\begin{equation}\label{20s21}
\tau v^2 -1=\frac{\xi_2} 2\left(\frac 1{(b^*)^2}+(b^*)^2\right) v.
\end{equation}
In this case we have a five diagonal matrix,
\begin{equation}\label{20s19}
v\sim S^2\Lambda_2+S\Lambda_1 +\Lambda_0 +\Lambda_1 S^{-1}
+\Lambda_2 S^{-2},
\end{equation}
but of a specific structure  (see Lemma 2.1).
Depending of ordering of
infinities
one of coefficients
$\lambda_0^{(2)},\lambda_1^{(2)}$ vanishes.
Say, $\lambda_1^{(2)}=0$, respectively
$\lambda_0^{(2)}\not=0$. Using $(b^*)^2\sim S^2$ we get
from
\eqref{20s21}
\begin{equation}\label{20s23}
\begin{split}
&\tau S^2\Lambda_2 S^2\Lambda_2=\frac{\xi_2}2 S^4\Lambda_2,\\
&\tau(S^2\Lambda_2 S\Lambda_1+S\Lambda_1 S^2\Lambda_2)=
\frac{\xi_2}2 S^3\Lambda_1,\\
&\tau(S^2\Lambda_2 \Lambda_0+S\Lambda_1 S\Lambda_1
+\Lambda_0 S^2\Lambda_2)=
\frac{\xi_2}2 S^2\Lambda_0,\\
&\tau(S^2\Lambda_2 \Lambda_1 S^{-1}+S\Lambda_1 \Lambda_0+
\Lambda_0 S \Lambda_1
+\Lambda_1 S\Lambda_2)=
\frac{\xi_2}2 S^2\Lambda_1 S^{-1},\\
&\tau(S^2\Lambda_2^2 S^{-2}+S\Lambda_1^2 S^{-1}+
\Lambda_0^2
+\Lambda_1^2+ \Lambda_2^2)-I=
\frac{\xi_2}2 (S^2\Lambda_2 S^{-2}+\Lambda_2).\\
\end{split}
\end{equation}
That is $\tau\lambda_0^{(2)}=\xi_2/2$ and the second relation
in \eqref{20s23} is an identity. Further,
\begin{equation}\label{20s24}
\tau\lambda_0^{(1)}\lambda_1^{(1)}
=-\frac{\xi}2\lambda_0^{(0)}=\frac{\xi}2\lambda_1^{(0)}
\end{equation}
and the fourth relation is an identity. Finally, from the last
equation we get
\begin{equation}\label{20s25}
(\lambda_0^{(0)})^2+
(\lambda_0^{(1)})^2+(\lambda_1^{(1)})^2
=1/\tau.
\end{equation}
Thus \eqref{20s24}, \eqref{20s25} are counterparts of
\eqref{20s20} in this case.

\subsection{The Renormalization Equation in the general case} For
$\pi:\overline{\CC}\to
\overline{\CC}$ given by $T_2$ (that is
$X_c\simeq\overline{\CC}$ and $Y_c\simeq\overline{\CC}$)
let $E$ be a system of intervals, a subset of
$[-\xi_2,\xi_2]$. As usual $F=\pi^{-1}(E)$. For a system of
points on $X_c\setminus E$
\begin{equation}\label{29s33}
P_1=\infty, \, P_2,\, ...,\, P_l,
\end{equation}
we define
$$
b_*^l=b_{P_1}...b_{P_l}.
$$

We consider $z_*$ given by the identical
map $X_c\to \CC$. The finite difference operator $\tilde J$
related to this symbol, a character $\alpha\in\Gamma^*_X$
and the ordering system of infinities \eqref{29s33}
is of the form
\begin{equation}\label{29s34}
\tilde J=
\begin{bmatrix}
&*&*&\hdots&* & & & & &\\
&*&*&\hdots&* & & & & &\\
& &*&\hdots&* & & & & &\\
& & &\ddots&\vdots & & & & &\\
& & & & *&*&*& \hdots&*  & \\
& & & & &*&*& \hdots&*  & \\
& & & & & &*& \hdots&*  & \\
& & & & & & &\ddots &\vdots &  \\
& & & & & & & &* &
\end{bmatrix},
\end{equation}
because, actually, $z_*(P_j)=\infty$ only for $j=1$.

The system of infinities on $Y_c$ is given by
\begin{equation}\label{29s35}
Q^{(1)}_1=\infty, \, Q^{(1)}_2=0,\, ...,\, Q^{(l)}_1,\, Q^{(l)}_2,
\end{equation}
where $\{Q_1^{(j)}, Q^{(j)}_2\}=\pi^{-1}(P_j)$, $j=1,...,l$.
Respectively,
$$
z^*=z_*\circ\pi=\tau v-\frac 1 v,\quad
(b^*)^2=b_*\circ\pi.
$$
We are interested to find a relation between
the initial operator $\tilde J$ and
the finite difference
operator $J$
with the symbol $v$
related to the character $\pi^*\alpha\in\Gamma_Y^*$ and
the ordering system of infinities \eqref{29s35}.
 Let us point out that
the symbol of $J$ is $v$, not $z^*$. It gives us an opportunity
to iterate the procedure, which appears to be a certain
rescaling.

Let us prove the following lemma.

\begin{lemma} Let
$T(z)=\frac{Q(z)}{P(z)}$. Define
\begin{equation*}
(\mathfrak L f)(x)=\frac 1 d \sum_{T(y)=x} f(y).
\end{equation*}
Then
\begin{equation*}
\mathfrak L \frac{1}{z-x}=\frac 1 d\frac{P'(z)}{P(z)}
+\frac 1 d\frac{T'(z)}{T(z)-x}.
\end{equation*}
\end{lemma}
\begin{proof}
Note that
\begin{equation*}
 \sum_{Q(y)-xP(y)=0} \frac{1}{z-y}
 =\frac{Q'(z)-xP'(z)}{Q(z)-xP(z)}
\end{equation*}
and then collect corresponding terms
\begin{equation*}
\begin{split}
 \frac{Q'(z)-xP'(z)}{Q(z)-xP(z)}=&
 \frac{Q'(z)-Q(z)P'(z)/P(z)+
 \{Q(z)-xP(z)\}P'(z)/P(z)}{Q(z)-xP(z)}\\
 =&\frac{P'(z)}{P(z)}+
  \frac{Q'(z)P(z)-Q(z)P'(z)}{P^2(z)}
 \frac{1}{T(z)-x}.
\end{split}
\end{equation*}
\end{proof}

\begin{theorem} The following relation
holds
\begin{equation}\label{re}
\mathbf{V}^*(z-J)^{-1}\mathbf{V}=
\frac 1 d\frac{P'(z)}{P(z)}
+\frac 1 d {T'(z)}(T(z)-\tilde J)^{-1},
\end{equation}
where $\tilde J=\mathbf{z}_*(\alpha)$
and  $J=\mathbf{v}(\pi^*\alpha)$.
\end{theorem}

\begin{proof}
We use the previous lemma and the functional representation
of all operators involved in \eqref{re}.
\end{proof}

\subsection{A vector representation of $H^2(\pi^*\alpha)$}
Due to
$$
\int_F|f|^2\omega^*=
\int_E\frac 1 2\sum_{\pi(Q)=P}|f|^2\omega_*
$$
we have a certain representation of $f\in H^2(\pi^*\alpha)$
\begin{equation}\label{30s36}
f\mapsto\frac 1{\sqrt{2}}
\begin{bmatrix} f(Q_1(P)) \\ f(Q_2(P))
\end{bmatrix}\mapsto
\begin{bmatrix} g_1(P)\\ g_2(P)
\end{bmatrix},
\end{equation}
where
\begin{equation*}
\begin{bmatrix} g_1(P)\\ g_2(P)
\end{bmatrix}=
\frac 1 2
\begin{bmatrix} f(Q_1(P)) +f(Q_2(P))\\ f(Q_1(P))-f(Q_2(P))
\end{bmatrix},
\end{equation*}
as a two--dimensional vector--function.
Let us describe analytical properties of this object
(metrically it is of $L^2$ with respect to $\omega_*$,
moreover the transformation is norm--preserved).

It is evident that the function
$g_1$
belongs to $H^2(\alpha)$. Consider the second function.
Let $c_{\pm}$ be critical points
$$
\pi'(c_{\pm})=0\Rightarrow
c_{\pm}=\pm\frac i{\sqrt{\tau}}.
$$
For a small circle $\gamma$ around the point
$\pi(c_{\pm})=\pm 2i\sqrt{\tau}$ we have
$g_2\circ\gamma=-g_2$. Let us introduce
$$
\Delta^2:= b_{\pi(c_+)} b_{\pi(c_-)}.
$$
Note that $\Delta\circ\gamma=-\Delta$. Further, since
the ratio
$$
\frac{b_{\pi(c_+)}} {b_{\pi(c_-)}}(z)=e^{iC}\frac
{z-\pi(c_+)} {z-\pi(c_-)}
$$
is singlevalued in $X_c\setminus E$ the both Green functions
have the same character of automorphy, which we denote by
$\nu\in \Gamma^*_X$. Thus, we get
\begin{equation}\label{12o38}
g_2=\Delta \tilde g_2,\
\text{where}\ \tilde g_2\in H^2(\alpha\nu^{-1}).
\end{equation}
We can summarize the result of this subsection as
\begin{theorem}
The transformation $f\mapsto g_1\oplus g_2$ given by
\eqref{30s36} is a unitary map from $H^2(\pi^*\alpha)$
to $H^2(\alpha)\oplus\Delta H^2(\alpha\nu^{-1})$.
Moreover,
$$
Vf\mapsto f\oplus 0,\quad f\in H^2(\alpha),
$$
where the isometry $V:H^2(\alpha)\to H^2(\pi^*\alpha)$
is defined by \eqref{24s19}. Also,
\begin{equation}\label{30s37}
k^{\pi^*\alpha}_{Q_{\pm}}=
k^{\alpha}_{P}\oplus
(\pm \Delta\overline{\Delta(P)}k^{\alpha\nu^{-1}}_{P})
\end{equation}
for the reproducing kernel
$k^{\pi^*\alpha}_{Q_{\pm}}\in H^2(\pi^*\alpha)$,
$\pi(Q_{\pm})=P$.
\end{theorem}

Extending this vector representation onto $L^2$--spaces we
get immediately
\begin{theorem}The transformation $f\mapsto g_1\oplus \tilde g_2$ given
by
\eqref{30s36} and \eqref{12o38} is a unitary map from $L^2(\pi^*\alpha)$
to $L^2(\alpha)\oplus L^2(\alpha\nu^{-1})$. With respect to this
representation the multiplication operator by $z^*=z_*\circ\pi$
is of the form
\begin{equation}\label{12o40}
\mathbf z^*(\pi^*\alpha)\mapsto
\begin{bmatrix}\mathbf z_*(\alpha)& 0\\
0&\mathbf z_*(\alpha\nu^{-1})
\end{bmatrix}.
\end{equation}
\end{theorem}

To get a similar representation for the multiplication operator by $v$
we need to introduce the following notations. Let us note that
$$
b_\infty^2\frac{z^2+4\tau}{\Delta^2}
$$
is an outer function in the domain $\overline{\CC}\setminus E$.
So, the square root of this function is well defined. We put
$$
b_\infty\phi:=\sqrt{
b_\infty^2\frac{z^2+4\tau}{\Delta^2}}.
$$
Since on the boundary of the domain
$$
\phi^2\Delta^2=z^2+4\tau=|\phi|^2
$$
(the second expression is positive on $E$)
we have
\begin{equation}\label{12o41}
\phi\Delta^2=\overline{\phi} \quad {\rm on}\ E.
\end{equation}

\begin{lemma}
The multiplication operator
$\phi: L^2(\alpha)\to L^2(\alpha\nu^{-1})$
with respect basises systems
related to the infinities
$\{P_1,...,P_l\}$ has as many diagonals as
$\mathbf z_*(\alpha)$ and $\mathbf z_*(\alpha\nu^{-1})$.
Moreover,
\begin{equation}\label{12o42}
\begin{split}
&\boldsymbol{\phi} \mathbf z_*(\alpha)=
\mathbf z_*(\alpha\nu^{-1}) \bphi,\\
&    \bphi^*\bphi   =
\mathbf z^2_*(\alpha)+4\tau,\\
&   \bphi \bphi^*    =
\mathbf z^2_*(\alpha\nu^{-1})+4\tau.
\end{split}
\end{equation}
\end{lemma}

\begin{proof} First of all $\phi$ is a character--automrphic function
with the character $\nu^{-1}$, therefore the multiplication
operator acts from  $L^2(\alpha)$ to $L^2(\alpha\nu^{-1})$.
Since $b_\infty\phi$ is an outer function, $\phi$ has a unique
pole at infinity, and, therefore, the operator $\bphi$
has the same structure over diagonal as the operator multiplication by
$z$. The adjoint operator has the symbol $\overline \phi$.
According to \eqref{12o41} it has analytic continuation
from the boundary inside the domain with the only pole at
infinity. Thus $\bphi^*$ is also of the same structure
over the main diagonal as
$\mathbf z_*(\alpha)$ or $\mathbf z_*(\alpha\nu^{-1})$.
Combining these two facts we get that the whole structure
of $\bphi$ coincide with the structure of the matrix
$\mathbf z_*(\alpha)$. Then, just comparing symbols of operators
on the left and right parts of \eqref{12o42}, we prove these identities.
\end{proof}

\begin{theorem} With respect to the
decomposition
$L^2(\pi^*\alpha)\simeq L^2(\alpha)\oplus L^2(\alpha\nu^{-1})$
the multiplication operator by $v$
is of the form
\begin{equation}\label{12o43}
\mathbf v(\pi^*\alpha)\simeq
\frac 1{2\tau}
\begin{bmatrix}\mathbf z_*(\alpha)& \bphi^*\\
\bphi &\mathbf z_*(\alpha\nu^{-1})
\end{bmatrix}.
\end{equation}
\end{theorem}

Let us mention that according to \eqref{12o42},
the operators given in \eqref{12o40} and \eqref{12o43}
commute and satisfy the identity, which is generated by the symbols
identity $z^*=\tau v-1/v$.


\subsection{
One sided matrices
and
the expanding transform $\pi(v)=\tau v-\frac{\tau-1}{v}$}

Note that in this normalization $v=1$ is the positive fixed point,
$\pi(1)=1$. Put
$E_0=[-1,1]$. For a continuous function $f$ on
$$
E_1=\pi^{-1}([-1,1])=[-1,-1+\frac{1}{\tau}]
\cup[1-\frac 1{\tau},1]
$$
we define
\begin{equation}
    (\cL f)(x)=\frac 1 2\sum_{\pi(v)=x}f(v).
\end{equation}
The conjugate operator acts on measures $$\cL^*:C(E_0)^*\to C(E_1)^*.$$

Let $f_0, f_1, f_2, ...$ be a certain orthonormal system with respect to a (positive) measure
$\nu\in C(E_0)^*$ then $$f_0\circ\pi, f_1\circ\pi, f_2\circ\pi, ...$$ is orthonormal system with respect
to $\mu:=\cL^*\nu$. Note that if the first system form basis in $L^2_{d\nu}$ the second one form basis in the set of "even" functions from $L^2_{d\mu}$, the functions that are invariant with respect to the substitution
$ f(-\frac{\tau-1}{\tau v})=f(v)$.

\begin{proof}[Example] Let $f_0, f_1, f_2, ...$ be orthonormal
polynomials in
$L^2_{d\nu}$.
$f_0\circ\pi, f_1\circ\pi, f_2\circ\pi, ...$ is a certain orthonormal system in $L^2_{d\mu}$ consisting of "polynomials" of $v$ and $1/v$, similarly to the systems that generate CMV matrices:
$$
    1, v,1/v, v^2, 1/v^2...
$$
Making a small modification in this procedure, we orthogonalize
$$
    1, \tau v+\frac{\tau-1}{v},  \tau v-\frac{\tau-1}{v},
    (\tau v)^2-\left(\frac{\tau-1}{v}\right)^2,
        (\tau v)^2+\left(\frac{\tau-1}{v}\right)^2\dots
$$
and denote the orthonormal system by $e_0, e_1, e_2,...$.

It is evident that
$$
e_{2k}=f_k\circ\pi
$$
 and
$$
e_{2k+1}(v)=\left(\tau v+\frac{\tau-1}{v}\right) g_k(\pi(v)),
$$
 where $g_k$ is also orthonormal system of
{\it polynomials} but with respect to the measure
$(x^2+4\tau(\tau-1))d\nu(x),
$ since
$$
\left(\tau v+\frac{\tau-1}{v}\right)^2=x^2+4\tau(\tau-1),\quad
{ for}\ x=\tau v-\frac{\tau-1}{v}.
$$
Let $J$ be the Jacobi matrix, corresponding to the multiplication operator
in $L^2_{d\nu}$ with respect to the basis of the orthonormal polynomials.

The given $J$ we want to describe the multiplication operator in $L^2_{d\mu}$ with respect to $\{e_k\}$.

We decompose
$L^2_{d\mu}$ onto even and odd subspaces. Then
$$
    \tau v-\frac{\tau-1}{v}\mapsto
    \left[
    \begin{array}{cc}
    J&0\\0 & J_*
    \end{array}
    \right],
$$
where $J_*$ is the Jacobi matrix corresponding to the measure $(x^2+4\tau(\tau-1))d\nu(x)$.

Further,
 $$
    \tau v+\frac{\tau-1}{v}\mapsto\left[\begin{array}{cc}
    0&\Phi^*\\ \Phi & 0
    \end{array}\right].
$$
It is quite evident that $\Phi$ is an upper triangular matrix.

We get that
$$
    v\mapsto\frac 1{2\tau}\left[\begin{array}{cc}
    J&\Phi^*\\ \Phi & J_*
    \end{array}\right],
$$
and
$$
    -1/v\mapsto\frac 1 {2(\tau-1)}\left[\begin{array}{cc}
    J&-\Phi^*\\ -\Phi & J_*
    \end{array}\right].
$$
Therefore
$$
    \left[\begin{array}{cc}
    J^2-\Phi^*\Phi &-J\Phi^*+\Phi^* J_*\\ \Phi J- J_*\Phi & {J_*}^2-\Phi\Phi^*
    \end{array}\right]=\left[\begin{array}{cc}
    -4\tau(\tau-1)& 0\\ 0 & -4\tau(\tau-1)
    \end{array}\right].
    $$
    Thus $\Phi$ can be found in the upper--lower triangular factorization
    \begin{equation}
    \Phi^*\Phi=J^2+4\tau(\tau-1),
\end{equation}
and for $J_*$ we have $J_* =\Phi J\Phi^{-1}$.

Thus for
$$
J=\left[
\begin{array}{cccc}
0&p_1& &\\
p_1&0& p_2& \\
 &\ddots&\ddots&\ddots
\end{array}\right]
$$
we have
$$
\pi^*(J)=\frac{1}{2\tau}\left[
\begin{array}{cccccccc}
0& & & & & & &\\
\lambda_0 &0 & & & & & & \\
 p_1&0 &0& & & & & \\
 0&p_1\frac{\lambda_1}{\lambda_0} &\lambda_1 &0 & & & & \\
 0& 0& p_2 &0 &0 & & & \\
 0& 0&0 &p_2\frac{\lambda_2}{\lambda_1}& \lambda_2& 0& &\\
 0& \frac{p_1p_2}{\lambda_0}&0 &0 &p_3 &0 &0 &\\
 & & \ddots& & &\ddots & \ddots&\\
\end{array}\right].
$$
The matrix is selfadjoint and $\lambda_n$ are defined by the recursion
\begin{equation}
\lambda_n^2= 4\tau (\tau-1) + p_{n+1}^2+p_n^2
 -\frac{p_n^2 p_{n-1}^2}{\lambda^2_{n-2}}.
\end{equation}
with the initial data
\begin{equation}
    \lambda_{0}^2=p_1^2+4\tau(\tau-1),\quad
    \lambda_1^2=p_1^2+ p_{2}^2+4\tau(\tau-1).
\end{equation}
\end{proof}

\begin{theorem}
Let $\nu$ be the spectral measure of $A$,
\begin{equation}\label{22def2}
    \int \frac {d\nu(x)}{x-z}=\langle 0|(A-z)^{-1}|0\rangle.
\end{equation}
 Then $\pi^*(A)$ is a self--adjoint operator
with the cyclic vector $|0\rangle$ and the spectral measure $\mu=\cL^* \nu$.
\end{theorem}

\begin{proof}
By  Definition \ref{22def} and \eqref{22def2} we have
\begin{equation}\label{223}
\begin{split}
\langle 0|(\pi^*(A)-z)^{-1}|0\rangle=& 2\tau\left\langle
\begin{bmatrix}A-2\tau z&\Phi^*\\
\Phi& A_*-2\tau z \end{bmatrix}^{-1}
\begin{bmatrix}|0\rangle\\0 \end{bmatrix},
\begin{bmatrix}|0\rangle\\0 \end{bmatrix}
\right\rangle\\
=&2\tau\langle 0|(A-2\tau z-\Phi^*(A_*-2\tau z)^{-1}\Phi)^{-1}|0\rangle\\
=&2\tau\langle 0|(A-2\tau z-\Phi^*\Phi(A-2\tau z)^{-1})^{-1}|0\rangle\\
=&2\tau\langle 0|(A-2\tau z-(A^2+4\tau(1-\tau))(A-2\tau z)^{-1})^{-1}|0\rangle\\
=&\int\frac{2\tau d\nu(x)}{x-2\tau z-(x^2+4\tau(1-\tau))(x-2\tau z)^{-1}}\\
=&\int\frac{(x-2\tau z)d\nu(x)}{2\tau z^2-2xz-2(1-\tau)}.
\end{split}
\end{equation}
Since
\begin{equation*}
    \left(\cL\frac 1{v-z}\right)(x)=\frac 1 2\sum_{\tau v-\frac{\tau-1}{v}=x}\frac 1{v-z}
    =\frac{x-2\tau z}{2\tau z^2-2xz-2(1-\tau)},
\end{equation*}
we get
\begin{equation*}
    \langle 0|(\pi^*(A)-z)^{-1}|0\rangle=\int\left(\cL
    \frac
    1{v-z}\right)(x)d\nu(x)=\int \frac 1{v-z}d(\cL^*\nu)(v)
\end{equation*}
and the theorem is proved.

\end{proof}

Using Ruelle's Theorem with respect to the map $\pi(v)$ we can
summarize our considerations with the following theorem.

\begin{theorem}\label{t27n}
The iterative procedure
$$
A_{n+1}=\pi^*(A_n)
$$
converges to the operator $A=\lim_{n\to \infty} A_n$ with the spectral measure $\mu$ which is
the eigen--measure for the Ruelle operator $\cL^*\mu=\mu$. The operator $A$ is the multiplication operator by independent variable in $L^2_{d\mu}$
with respect to the following basis
$$
e_{2k}(v)=e_k(\pi(v)),
$$
$$
e_{2k+1}(v)=\left(\tau v+\frac{\tau-1}{v}\right)\sum_{j=0}^k c^k_j e_j(\pi(v)), \quad e_0(v)=1,
$$
where the coefficients $c^k_j$ with $c^k_k>0$ are uniquely determined due to the orthogonality
condition
$\langle e_{2k+1}, e_l\rangle=\delta_{2k+1,l}$, $l\le 2k+1$.  Moreover,
$e_k(v)$ is a rational function of $v$ such that
$
e_k(A)|0\rangle=|k\rangle,
$
and
\begin{equation}
\left[\begin{array}{cccc}
c^0_0&c^1_0&c_0^2 &\dots\\
0&c^1_1& c^2_1& \dots\\
 &\ddots&\ddots&\ddots
\end{array}\right]=\Phi^{-1}.
\end{equation}
\end{theorem}

\bibliographystyle{amsplain}

Addresses:\\[.1cm]
\noindent
Franz Peherstorfer \\
Abteilung f\"ur Dynamische Systeme \\
und Approximationstheorie\\
Universit\"at Linz\\
4040 Linz, Austria\\
{Franz.Peherstorfer@jku.at}\\

\noindent
Peter Yuditskii\\
Abteilung f\"ur Dynamische Systeme \\
und Approximationstheorie\\
Universit\"at Linz\\
4040 Linz, Austria\\
 Petro.Yudytskiy@jku.at

\end{document}